\documentclass[a4paper,12pt]{article}

\advance\textheight2.9cm        
\advance\topmargin-2.0cm
\advance\textwidth2.6cm
\advance\evensidemargin-1.6cm

\usepackage{amsfonts} 
\usepackage{graphicx}
\usepackage{subcaption}
\usepackage{url}
\usepackage{cite}

\usepackage{theorem}
\usepackage{amsfonts}
\usepackage{amssymb} 
\usepackage{amsmath}
\usepackage{longtable}
\usepackage{array}
\usepackage{bm}  
\usepackage{comment}
\def\comment#1{{}}	

\usepackage{theorem}

\newtheorem{lem}{Lemma}

\newtheorem{conj}{Conjecture}

\def\Box{{\hbox{\raisebox{0.0em}{\rlap{$\sqcap$}}\kern0em%
            \raisebox{-0.0em}{$\sqcup$}}} } 

\newenvironment{proof}{{\it Proof. }}{\nopagebreak\hspace*{0.5cm}\hfill$\Box$\vspace{0.5cm}}

\def\bepf{\begin{proof}}
\def\epf{\end{proof}}

\def\eeq{\end{equation}}
\def\lbeq#1{\begin{equation} \label{#1}}
\def\bary{\begin{array}}
\def\eary{\end{array}}

\def\ealg{\end{alg}}

\usepackage{array,multirow,graphicx}
\usepackage[utf8]{inputenc}

\usepackage{graphicx}

\usepackage{colortbl}

\usepackage[table]{xcolor}

\usepackage{fancyhdr}

\usepackage{hyperref}
\usepackage{cleveref}

\usepackage{array,multirow,graphicx}

\usepackage{etex}

\usepackage{longtable}
\usepackage{lscape}
\usepackage{hhline}
\usepackage{caption}
\usepackage{blindtext}
\usepackage{lmodern}
\usepackage{textcomp}
\usepackage{chngcntr}
\usepackage{tikz}
\usepackage{csquotes}
\usepackage[english]{babel}
\usepackage[utf8]{inputenc}
\usepackage{algorithm}
\usepackage{xcolor}
\usepackage{ulem}
\usepackage{cancel}
\usepackage{multicol}
\usepackage{algpseudocode}

\fancypagestyle{plain}{\fancyhf{}\fancyfoot[C]{\bfseries \thepage}}

\newcolumntype{?}{!{\vrule width 1pt}}

\usepackage{adjustbox}

\newcounter{subsubsubsection}[subsubsection]
\def\subsubsubsectionmark#1{}

\def\subsubsubsection{\@startsection
	{subsubsubsection}{4}{\z@} {-3.25ex plus -1
		ex minus -.2ex}{1.5ex plus .2ex}{\normalsize\bf}}
\def\l@subsubsubsection{\@dottedtocline{4}{4.8em}
	{4.2em}}

\newcolumntype{R}[2]{%
	>{\adjustbox{angle=#1,lap=\width-(#2)}\bgroup}%
	l%
	<{\egroup}%
}

\usetikzlibrary{matrix,shapes,arrows,positioning,chains}



\usepackage{comment}
\usepackage{tikz}
\usetikzlibrary{matrix,shapes,arrows,positioning,chains}

\definecolor{ao(english)}{rgb}{0.0, 0.5, 0.0}

 \numberwithin{equation}{section}
 
 \newcommand{\bd}{\mathbf{d}}
 \newcommand{\bx}{\mathbf{x}}

\newcommand{\bH}{\mathbf{H}}
 \newcommand{\bQ}{\mathbf{Q}}
 \newcommand{\bg}{\mathbf{g}}
 \newcommand{\by}{\mathbf{y}}
\newcommand{\bs}{\mathbf{s}}

\newcommand{\bP}{\mathbf{P}}

\begin{document}

\begin{center}

{\Large \bf From a Scalar to a Matrix Setting for the Dai--Liao Parameter}

\vspace{0.5cm}

{\large \bf Saman Babaie--Kafaki}
\centerline{\sl Faculty of Engineering, Free University of Bozen--Bolzano}
\centerline{\sl NOI Techpark, Via Bruno Buozzi 1, 39100 Bolzano (BZ), Italy}
\centerline{\sl email: saman.babaiekafaki@unibz.it}

\vspace{0.5cm}

{\large \bf Morteza Kimiaei}
\centerline{\sl Fakult\"at f\"ur Mathematik, Universit\"at Wien}
\centerline{\sl Oskar--Morgenstern--Platz 1, A--1090 Wien, Austria}
\centerline{\sl email: kimiaeim83@univie.ac.at}

\vspace{0.5cm}

{\large \bf Zohre Aminifard}
\centerline{\sl Université Catholique de Louvain (UCLouvain)}
\centerline{\sl Institute of Information and Communication Technologies, Electronics and} 
\centerline{\sl  Applied Mathematics, Place du Levant 3, B--1348, Louvain--la--Neuve, Belgium}
\centerline{\sl email: zohreh.aminifard@uclouvain.be}

\vspace{0.5cm}

\end{center}

{\bf Abstract.} As is well known, both the numerical performance and the theoretical properties of the Dai--Liao conjugate gradient algorithm are highly dependent on the adjustment of its key parameter. Here, we first employ the well-known Harmonic--Geometric--Arithmetic--Quadratic mean inequality to reinforce the optimality of two previously proposed scalar adaptive settings of the Dai--Liao parameter. In other words, we show how these parameter choices are capable of enhancing the well-conditioning of the Dai--Liao search direction matrix by shrinking the intervals containing its singular values. A similar analysis is also carried out for scaled memoryless quasi--Newton updating formulas in order to further justify the optimality of two classical scaling parameters associated with these updates. Then, as the main contribution of this work, we move from the classical scalar setting of the Dai--Liao parameter to a matrix setting aimed at enhancing flexibility and diversity within optimization methods. In particular, we show that, under such a matrix formulation of the parameter, a well-known three-term conjugate gradient algorithm emerges as a member of the proposed extended Dai--Liao class of algorithms while enjoying several computationally attractive properties. Among these, the well-conditioning of the associated search direction matrix is especially noteworthy, as well as the ability to make more explicit use of the second-order information of the model.  Finally, to provide practical evidence supporting the proposed matrix setting of the Dai--Liao parameter, we conduct a series of numerical experiments on standard benchmark test problems and report the results in detail. Generally speaking, the proposed framework is shown to retain both theoretical soundness and computational reliability. 



\vspace{0.2cm} {\em 2000 AMS Subject Classification: 65K05; 90C53; 15A23.}

\hfill \today

\tableofcontents

\section{Introduction}\label{Introduction}


The growing influence of machine learning and data-driven technologies has further elevated the importance of \textbf{Nonlinear Programming} (NLP) in modern scientific computing \cite{Elden}. This trend has motivated extensive research aimed at improving both the theoretical foundations and computational aspects of optimization methods. Although numerous advances have enriched the available NLP tools, emerging applications and increasingly sophisticated problem structures require methodologies that offer greater flexibility, adaptability, and robustness. Among these research directions, limited-memory NLP methods occupy a prominent position due to the widespread occurrence of large-scale problems across numerous application domains. Meanwhile, among NLP models, unconstrained optimization plays a fundamental role in practical applications, as many NLP methods repeatedly solve unconstrained optimization subproblems in order to update the underlying approximations  \cite{Andreibook}.

\subsection{Preliminaries}

As is well known, to solve the unconstrained optimization problem 
\begin{equation}\label{UP}
    \displaystyle\min_{\bx\in \mathbb{R}^n}\ f(\bx),
\end{equation}
\textbf{Conjugate Gradient} (\texttt{CG}) iterations are recursively updated by
\begin{equation}\label{xk}
 \bx_{k+1} := \bx_k + \bs_k,
\end{equation}
for all $k \in \mathbb{Z}^+ := \{0,1,2,\dots\}$, starting from an initial guess $\bx_0 \in \mathbb{R}^n$, where $\bs_k \in \mathbb{R}^n$ is called the $k$th step of the algorithm. In this context, $\bs_k$ is typically determined as a scalar-vector product of the form $\bs_k := \alpha_k \bd_k$, where
\begin{itemize}
\item[$\blacktriangleright$] $\bd_k \in \mathbb{R}^n$ is the (descent) search direction, usually defined by
\begin{equation}\label{dk}
 \bd_{0} := -\bg_{0}, \qquad \bd_{k+1} := -\bg_{k+1} + \beta_k \bd_k, \qquad \forall k \in \mathbb{Z}^+,
\end{equation}
with $\bg_k := \nabla f(\bx_k)$, and $\beta_k \in \mathbb{R}$ referred to as the \texttt{CG} (update) parameter;
\item[$\blacktriangleright$] $\alpha_k \in (0,+\infty)$ is the step size, typically computed via a \textbf{Line Search} (\texttt{LiS}) strategy along $\bd_k$ to locally ensure sufficient decrease in the cost function $f$ \cite{Andreibook}.
\end{itemize}

As their name suggests, \texttt{CG} algorithms were initially based on the so-called ``conjugacy condition". Although for quadratic versions of \eqref{UP} such a condition can be guaranteed for all the search directions (under the exact {\tt LiS}) \cite{Andreibook}, for general cost functions the conjugacy condition cannot be ensured. Some \texttt{CG} algorithms have their roots in a successive conjugacy condition. In particular, the \textbf{Hestenes--Stiefel} (\texttt{HS}) \cite{HS} parameter, i.e.,
\begin{equation*}
    \beta_k^{\texttt{HS}} := \dfrac{\bg_{k+1}^\top \by_k}{\bd_k^\top \by_k},
\end{equation*}
is derived from the condition
\begin{equation}\label{conjugacy0}
    \bd_{k+1}^\top \by_k = 0,
\end{equation}
with $\by_k := \bg_{k+1} - \bg_k$ as the gradient displacement vector \cite{SBKRairo}.  Although the \texttt{HS} method has been traditionally categorized among the classical \texttt{CG} algorithms with acceptable computational performance \cite{AndreiComparison}, its global convergence requires strong assumptions to be guaranteed \cite{GN}, especially since the method fails to ensure \textbf{Sufficient Descent Condition} (SDC).


At the beginning of the current century, the \textbf{Dai--Liao} (\texttt{DL}) algorithm was developed by extending \eqref{conjugacy0} \cite{DaiLiaoNCG}, benefiting also from \textbf{Quasi–Newton} (\texttt{QN}) aspects, in the sense that
\begin{equation}\label{conjugacyDL}
\bd_{k+1}^\top\by_k = -t\bg_{k+1}^\top\bs_k,
\end{equation}
in which the scalar $t>0$ is referred to as the \texttt{DL} parameter. Then, using \eqref{dk} and \eqref{conjugacyDL}, an extended version of $\beta_k^{\texttt{HS}}$ has been formulated as
\begin{equation}\label{betdl}
\beta_k^{\texttt{DL}} := \beta_k^{\texttt{HS}} - t\dfrac{\bg_{k+1}^\top\bs_k}{\bd_k^\top\by_k},
\end{equation}
which has the potential to ensure global convergence for uniformly (strongly) convex cost functions, as well as its restricted version given by
\begin{equation*}
\beta_k^{\texttt{DL}+} := \beta_k^{\texttt{HS}+} - t\dfrac{\bg_{k+1}^\top\bs_k}{\bd_k^\top\by_k},
\end{equation*}
with $\beta_k^{\texttt{HS}+} := \max\left\{\beta_k^{\texttt{HS}}, 0\right\}$, which has the potential to ensure global convergence for general cost functions \cite{DaiLiaoNCG}. A careful review of the literature reveals that the \texttt{DL} algorithm can be regarded as a landmark approach due to its remarkable flexibility in solving unconstrained optimization problems \cite{SBKRairo}. In particular, the method guarantees the SDC when
\begin{equation}\label{DLtau}
t \geq \varsigma\dfrac{\|\by_k\|^2}{\bs_k^\top\by_k},\qquad \text{for some real constant } \varsigma > \dfrac{1}{4},
\end{equation}
where $\|.\|$ signifies the $\ell_2$ norm \cite{HagerZhang}. Hence, an effective \texttt{LiS} strategy can further enhance the performance of the \texttt{DL} algorithm.

Nowadays, the Wolfe \texttt{LiS} conditions are widely adopted as an effective hybridization of the Armijo condition, i.e.,
\begin{equation}\label{wolfe1}
f(\bx_k + \alpha_k \bd_k) - f(\bx_k) \leq \rho_1 \alpha_k\left(\bg_k^\top \bd_k\right),\qquad \rho_1 \in (0,1),
\end{equation}
together with the curvature condition, i.e.,
\begin{equation}\label{wolfe2}
\nabla f(\bx_k + \alpha_k \bd_k)^\top \bd_k \geq \rho_2 \left(\bg_k^\top \bd_k\right),\qquad \rho_2 \in (\rho_1,1).
\end{equation}
As a consequence of inequality \eqref{wolfe2}, when $\bd_k$ is a descent direction, it follows that
\begin{equation}\label{sy}
    \bd_k^\top \by_k > 0,
\end{equation}
a key property that ensures the well-definedness of the \texttt{DL} parameter. In the subsequent analysis, we assume that the \texttt{LiS} strategy satisfies the Wolfe conditions \eqref{wolfe1} and \eqref{wolfe2}.

As is well known, both the \texttt{LiS} conditions and the SDC are key components in the global convergence analysis of \texttt{CG} algorithms \cite{CGConvDai}. Meanwhile, two common assumptions that are typically required to establish the convergence of \texttt{CG} algorithms are as follows \cite{SugikiNarushimaYabe}:
\begin{enumerate}
\item[$(i)$] boundedness of the level set $\mathfrak{L}_0 := \{\bx : f(\bx) \leq f(\bx_0)\}$;
\item[$(ii)$] Lipschitz continuity of $\nabla f$ in a neighborhood $\mathfrak{L}$ of the set $\mathfrak{L}_0$.
\end{enumerate}
These assumptions are considered to hold throughout this study. Also, the following version of Lemma~3.1 in \cite{SugikiNarushimaYabe} serves as the foundation for our forthcoming convergence analysis.

\begin{lem}\label{lemconv}
Consider an iterative algorithm of the form \eqref{xk}--\eqref{dk}, where the SDC holds, and the \texttt{LiS} is performed to satisfy the Wolfe conditions \eqref{wolfe1} and \eqref{wolfe2}. If
\begin{equation*}
    \displaystyle\sum_{k=0}^{+\infty}\dfrac{1}{\|\bd_k\|^2}=+\infty,
\end{equation*}
then the method converges globally in the sense that
\begin{equation*}
    \liminf_{k \to \infty} \| \bg_k \| = 0.
\end{equation*}
\end{lem}

As a result of Lemma~\ref{lemconv}, the uniform boundedness of the sequence $\left\{\|\bd_k\|\right\}_{k\geq0}$ can serve as a sufficient condition for global convergence. Meanwhile, if a \texttt{CG} algorithm fails to ensure such boundedness or cannot guarantee the SDC, a safeguard can be imposed on the algorithm in the form of a restart using the steepest descent direction \cite{Kimiaei2026}.

\subsection{Contribution of the Study} 
Here, before presenting the main contribution of this work, first by considering the relationship between the well-known mathematical means widely used in science and engineering \cite{Peressini}, i.e., \textbf{Arithmetic Mean} (AM), \textbf{Geometric Mean} (GM), \textbf{Harmonic Mean} (HM), and \textbf{Quadratic Mean} (QM) \cite{Inequalities}, we provide a further optimality feature for two previously proposed choices for the \texttt{DL} parameter $t$. More precisely, our analysis is fundamentally based on the ordered form of the mentioned means in the framework of HM--GM--AM--QM inequality \cite{Inequalities}, incorporating the singular values of the \texttt{DL}  \textbf{Search Direction Matrix} (SDM)  in a way to promote well-conditioning.  We also show that a similar analysis can be carried out for scaled memoryless \textbf{Quasi--Newton} (\texttt{QN}) updating formulas to further justify the optimality of two classical scaling parameters associated with these updates.

The main contribution of this work lies in extending the classical scalar setting of the \texttt{DL} parameter to a more general matrix framework. In particular, we introduce a dimensionality-expansion strategy within the \texttt{DL} algorithm by employing limited-memory matrix structures. This extension offers a flexible mechanism for exploiting additional structural information of the underlying model while simultaneously allowing the integration of fundamental analytical properties, such as conjugacy \cite{SBKRairo}. Moreover, it enables the establishment of potential connections between the \texttt{DL} method and other \texttt{CG} algorithms, such as the efficient three-term \texttt{CG} method proposed by \textbf{Zhang, Zhou, and Li} ({\tt ZZL}) \cite{ZhangZhouLiOMS}. From an accuracy standpoint, the proposed matrix setting for the \texttt{DL} parameter can be seamlessly incorporated into all versions of the method that utilize modified secant equations \cite{SBKRairo}, which exploit both cost function values and gradient information. In terms of consistency, all \texttt{CG} algorithms derived from the \texttt{DL} framework based on classical \texttt{CG} parameters can also benefit from the proposed matrix adjustment \cite{SBKRairo}.

\subsection{Organization of the Work} The organization of the work is as follows. In Section \ref{CG1}, we analyze the optimality of two previously proposed scalar settings of the \texttt{DL} parameter through a singular value analysis based on the HM--GM--AM--QM inequality. A similar analysis is also carried out for scaled memoryless \texttt{QN} updating formulas. Section \ref{MainCG} is devoted to the development of a matrix setting for the \texttt{DL} parameter. In particular, we employ a specific rank-one update of a scaled version of the identity matrix and investigate several attractive properties of the resulting parametric framework. Numerical experiments on standard benchmark test problems are presented in Section \ref{Numerical} to provide practical evidence supporting the proposed matrix setting of the \texttt{DL} parameter. Finally, in Section \ref{Conclusions}, we summarize the main findings of the study and outline several directions for future research.

\section{Parametric Setting for Memoryless Optimization Algorithms Based on Classical Inequalities}\label{CG1}

Mathematical inequalities have traditionally played significant roles in analyzing and modeling algorithmic behaviors across a wide range of applications in both science and engineering. In particular, from an algorithmic viewpoint, algebraic inequalities have been regarded as the foundation of complexity and error analysis, serving as an essential tool for assessing the stability of computational implementations as well \cite{Watkins}. Moreover, from a modeling perspective, inequalities often provide reliable means of transforming complex models into more tractable forms through various relaxation schemes \cite{SunYuan}. Consequently, it has now become a common practice to exploit classic mathematical inequalities for making judicious parametric settings.

In the field of optimization, inequalities have been extensively employed to establish effective lower or upper bounds for various applications. The effectiveness of such parametric bound inequalities is often evaluated by their proximity to the actual values of the corresponding parameters---a property commonly referred to as the ``sharpness" or ``pointedness" of an inequality \cite{Watkins}. In other words, the closer a bound is to the true parameter value, the more reliable the inequality becomes. However, the structural simplicity of a bound inequality is also of considerable importance, as it is directly related to the computational cost of algebraic manipulation. Therefore, when utilizing an inequality, it is essential to strike an appropriate balance between sharpness (as an indicator of accuracy) and formulaic simplicity (as an indicator of efficiency)---a challenge that can often be difficult to handle. Moreover, with the increasing prevalence of high-dimensional models in the modern era of big data, greater emphasis is now placed on the structural simplicity of bound inequalities \cite{SBKRairo}.

In real matrix spaces, the simplicity of inequalities is closely related to the choice of matrix norm. As is well known, the $\ell_2$ norm is theoretically one of the most popular norms; however, from a computational perspective, it is often considered costly because of being deeply connected to the singular values of the matrix  \cite{Watkins}. On the other hand, the well-known $\ell_1$ and $\ell_\infty$ matrix norms are numerically less expensive, but they do not define smooth mappings---a possible challenge from the viewpoint of algebraic manipulation \cite{SunYuan}. It is worth noting that all (induced) matrix norms are analytically equivalent, so it is possible to select a norm that is better aligned with specific objectives, providing a valuable degree of flexibility. To address the challenges associated with choosing an appropriate matrix norm for formulating an inequality, the Frobenius norm can sometimes be advantageous, as it is a direct extension of the Euclidean vector norm \cite{SunYuan}.

A careful review of the continuous optimization literature highlights the role of singular value and eigenvalue analyses in both convergence and the parameter tuning of algorithms \cite{SBKRairo}. In particular, even in limited-memory optimization algorithms, the structure of the SDM has a significant impact on algorithmic behavior. In such cases, it is important to incorporate as much information as possible from the SDM to enable more reliable parameter adjustments. From this perspective, the Frobenius norm is more advantageous than the other matrix norms, as it accounts for all the singular values (or eigenvalues, in the case of positive definiteness) of a matrix, unlike the Euclidean norm, which considers only the largest singular value \cite{Watkins}.

In what follows, the parameter adjustment of two well-known classes of memoryless optimization methods---specifically, \texttt{DL} and limited-memory \texttt{QN} algorithms---is investigated through the HM--GM--AM--QM inequality \cite{Inequalities}. Since singular values of the SDM of  \texttt{DL} method and eigenvalues of certain updating formulas of \texttt{QN} algorithms have already been characterized \cite{BabaieGhanbariEJOR,SBKJOTA3}, the proposed approach provides an effective plan to simultaneously and directly incorporate this information, thereby offering stronger justifications for the classical parameter tunings of these algorithms. From a conditioning perspective, the proposed framework focuses on the distribution of the singular values or eigenvalues of the SDM  \cite{Watkins}. Moreover, the approach emphasizes how such a formulation facilitates algebraic manipulations---such as differentiation---while maintaining sufficient flexibility.

\subsection{The  Mean Inequalities}

As an essential tool in basic statistics, representing the central tendency of a data set, the concept of the mean was originally defined in terms of AM. For a set of positive numbers $\mathcal{A}=\{a_1,a_2,\dots,a_n\}$, it is defined by
\begin{equation*}
    \textrm{AM}(\mathcal{A}) := \dfrac{a_1+a_2+\dots+ a_n}{n}=\dfrac{1}{n}\sum_{k=1}^na_k,
\end{equation*}
which represents the average of the given numbers and can be considered a computationally inexpensive formula. Based on the algebraic concept of the arithmetic operator, GM (defined by the operator ``$\times$") can be regarded as an extension of AM (defined by the operator ``+"); that is
\begin{equation*}
    \textrm{GM}(\mathcal{A}) := \left(a_1\times a_2\times\dots \times a_n\right)^{\frac{1}{n}}=\sqrt[n]{\prod_{k=1}^na_k},
\end{equation*}
often used in specific fields of economics and finance. Meanwhile, being useful especially in certain applications in physics, HM is defined by
\begin{equation*}
    \textrm{HM}(\mathcal{A}) := \dfrac{n}{\dfrac{1}{a_1} + \dfrac{1}{a_2}+\dots+ \dfrac{1}{a_n}}=\dfrac{1}{\textrm{AM}\left(\{a_k^{-1}\}_{k=1}^n\right)}.
\end{equation*}
In addition, formulated by
\begin{equation*}
    \textrm{QM}(\mathcal{A}) := \sqrt{\dfrac{a_1^2+a_2^2+\dots+ a_n^2}{n}},
\end{equation*}
QM is mainly used as a distribution evaluation factor in basic statistics (in terms of standard deviation), and also, as a measure to perform least-squares estimations \cite{Arashi}.

As already mentioned, inequalities play a crucial role in various branches of mathematics, particularly in optimization, where several fundamental inequalities stem from the concept of convexity. It is worth noting that the relationship between AM and GM in the context of the AM--GM inequality has received considerable attention in optimization. In particular, from a theoretical perspective, it can be regarded as the foundation of geometric programming \cite{Peressini}. On the other hand, the AM--GM inequality has been directly employed by Dennis and Wolkowicz \cite{DennisWolkowicz} to define a measure for analyzing the formulaic optimality as well as the convergence features of the \texttt{QN} algorithms. More recently, their approach has been extended to propose another measure function aimed at improving the condition number of the elements of nonnegative matrix factorization \cite{DargahiMMA}, by exploiting the relationship among the AM, GM, and HM. Here, the following comprehensive result is of vital importance, as it presents the way in which the aforementioned means can be ordered according to their magnitudes.

\begin{lem}\cite{Inequalities} Let $\mathcal{A}=\{a_1,a_2,\dots,a_n\}$ be a set of $n$ positive numbers. Then
\begin{equation}\label{HMGMAMQM}
    \textrm{HM}\left(\mathcal{A}\right)\leq \textrm{GM}\left(\mathcal{A}\right)\leq \textrm{AM}\left(\mathcal{A}\right)\leq \textrm{QM}\left(\mathcal{A}\right).
\end{equation}
Moreover, in (\ref{HMGMAMQM}), equality holds if and only if $a_1=a_2=\dots=a_n$.
\end{lem}

As can be seen, when $n=2$, inequality (\ref{HMGMAMQM}) reduces to
\begin{equation}\label{HMGMAMQM2}
\dfrac{2}{\dfrac{1}{a_1} + \dfrac{1}{a_2}} \leq \sqrt{a_1a_2} \leq \dfrac{a_1 + a_2}{2} \leq \sqrt{\dfrac{a_1^2 + a_2^2}{2}},
\end{equation}
where equality holds when $a_1 = a_2$. In the following, the role of (\ref{HMGMAMQM2}) in the parameter setting of two classes of memoryless optimization algorithms is highlighted. In general, here we consider the \texttt{LiS}  algorithms for solving (\ref{UP}) with the search directions which can be formulated as 
\begin{equation}\label{dk11}
\bd_0= \bg_0,\qquad \bd_{k+1}=-\bH_{k+1} \bg_{k+1}, \qquad \forall k \in \mathbb{Z}^+,
\end{equation}
where $\bH_k\in\mathbb{R}^{n\times n}$ is the SDM. The performance of such algorithms is highly dependent on the structure of their SDM. As observed, $\bg_{k+1}$ contains the first-order information of the cost function, while $\bH_{k+1}$ is often designed to explicitly (as in \texttt{QN} algorithms) or implicitly (as in \texttt{CG} algorithms, like \texttt{DL}) incorporate the second-order information of the model \cite{SunYuan}. Moreover, modified secant equations are quite helpful in engaging the available cost function values as well \cite{SBKRairo}.

As a fundamental factor that may affect the performance of an algorithm defined by a matrix-based iterative formula like (\ref{dk11}), the magnitude of the condition number of $\bH_{k+1}$ is a significant subject of investigation \cite{Watkins}. In particular, due to the frequent occurrence of large-scale models in real-world applications, devising limited-memory versions of $\bd_{k+1}$ has become crucial. This issue has often been addressed by formulating $\bH_{k+1}$ as a rank-one or rank-two update of some scaled forms of the identity matrix \cite{SBKRairo}. For such so-called vector-structured SDMs, in addition to computing $\bd_{k+1}$ via a few vector inner products, it is also possible to perform singular value or eigenvalue analyses \cite{SBKJOTA3,BabaieGhanbariEJOR}. Consequently, the distribution pattern of their singular values or eigenvalues can be optimized to make the corresponding SDMs as well-conditioned as possible.

\subsection{Optimal Scalar Parametric Settings for the Dai--Liao Algorithm}\label{DL}


As established in the literature, for the \texttt{DL} algorithm, the SDM is defined by
\begin{equation*}
\bP_{k+1}:=\mathbf{I}-\dfrac{\bs_k\by_k^\top}{\bs_k^\top\by_k} + t\dfrac{\bs_k\bs_k^\top}{\bs_k^\top\by_k}, \qquad \forall k \in \mathbb{Z}^+,
\end{equation*}
starting with $\bP_{0}=\mathbf{I}$, in the sense that, from  \eqref{dk} and \eqref{betdl}, the \texttt{DL} search direction can be written as $\bd_{k+1}=-\bP_{k+1}\bg_{k+1}$. To the best of our knowledge, determining an optimal formula for adjusting the \texttt{DL} parameter remains a well-known open problem in the literature of \texttt{CG} algorithms \cite{AndreiOpenProb}, which has attracted significant scholarly interest from both theoretical and computational perspectives \cite{SBKRairo}. In the following, we investigate the singular value distribution of $\bP_{k+1}$ using (\ref{HMGMAMQM2}).

As shown in \cite{BabaieGhanbariEJOR}, $\bP_{k+1}$ possesses $n-2$ unit singular values, along with two non-unit ones, denoted by $\sigma_k^-$ and $\sigma_k^+$, which are given by
\begin{eqnarray*}
 \sigma_k^\pm &:=& \sqrt{\dfrac{t^2\|\bs_k\|^4 + 2t(\bs_k^\top \by_k)\|\bs_k\|^2 + \|\bs_k\|^2\|\by_k\|^2}{4(\bs_k^\top \by_k)^2}} \\
      &&\qquad\qquad\pm\quad
        \sqrt{\dfrac{t^2\|\bs_k\|^4 - 2t(\bs_k^\top \by_k)\|\bs_k\|^2 + \|\bs_k\|^2\|\by_k\|^2}{4(\bs_k^\top \by_k)^2}}.
\end{eqnarray*}
It is worth noting that the Dennis--Wolkowicz measure function is actually defined as the quotient of AM to GM of the eigenvalues of a positive definite matrix (which are also its singular values) \cite{DennisWolkowicz}. Hence, it can be regarded as a general framework for well-conditioning, aiming to minimize the difference between the larger and the smaller mean of the singular values of a given matrix.

Here, note that for $\bP_{k+1}$, since
\begin{equation}\label{sigma1<1}
    \sigma_k^- \leq 1 \leq \sigma_k^+,
\end{equation}
it is reasonable to focus on reducing the gap between these two limiting singular values to improve the condition number of $\bP_{k+1}$. In this context, by defining $\mathcal{S}_k := \{\sigma_k^-, \sigma_k^+\}$, we can write
\begin{eqnarray*}
   \textrm{HM}(\mathcal{S}_k)&=& \dfrac{2t\|\bs_k\|^2}{\sqrt{t^2\|\bs_k\|^4 + 2t(\bs_k^\top \by_k)\|\bs_k\|^2 + \|\bs_k\|^2\|\by_k\|^2}}, \\
   \textrm{GM}(\mathcal{S}_k)&=& \sqrt{t\dfrac{\|\bs_k\|^2}{\bs_k^\top \by_k}}, \\
   \textrm{AM}(\mathcal{S}_k)&=& \sqrt{\dfrac{t^2\|\bs_k\|^4 + 2t(\bs_k^\top \by_k)\|\bs_k\|^2 + \|\bs_k\|^2\|\by_k\|^2}{4(\bs_k^\top \by_k)^2}}, \\
   \textrm{QM}(\mathcal{S}_k)&=& \sqrt{\dfrac{t^2\|\bs_k\|^4 + \|\bs_k\|^2\|\by_k\|^2}{2(\bs_k^\top \by_k)^2}}.
\end{eqnarray*}

Now, in a fractional programming framework similar to that of \cite{DennisWolkowicz}, from (\ref{HMGMAMQM2}) we have
\begin{eqnarray*}
   t_k^* := \dfrac{\|\by_k\|}{\|\bs_k\|} &=& \arg\min_{t} \dfrac{\textrm{AM}(\mathcal{S}_k)}{\textrm{GM}(\mathcal{S}_k)}
   = \arg\min_{t} \dfrac{\textrm{GM}(\mathcal{S}_k)}{\textrm{HM}(\mathcal{S}_k)} \\
   & =& \arg\min_{t} \dfrac{\textrm{AM}(\mathcal{S}_k)}{\textrm{HM}(\mathcal{S}_k)}
   = \arg\min_{t} \dfrac{\textrm{QM}(\mathcal{S}_k)}{\textrm{GM}(\mathcal{S}_k)}\\
   &=& \arg\min_{t} \dfrac{\textrm{QM}(\mathcal{S}_k)}{\textrm{HM}(\mathcal{S}_k)}= \arg\min_{t} \dfrac{\textrm{QM}(\mathcal{S}_k)}{\textrm{AM}(\mathcal{S}_k)}.
\end{eqnarray*}
In particular, regarding the last equality, it should be noted that
\begin{equation*}
    \dfrac{\textrm{AM}(\mathcal{S}_k)}{\textrm{QM}(\mathcal{S}_k)}
    = \sqrt{\dfrac{1}{2} + \dfrac{t(\bs_k^\top \by_k)}{t^2\|\bs_k\|^2 + \|\by_k\|^2}},
\end{equation*}
which attains its maximum at $t_k^*$.

It is worth noting that $t_k^*$ was also obtained in \cite{BabaieGhanbariEJOR} as the unique minimizer of an upper bound of the spectral condition number of $\bP_{k+1}$. Although it does not exceed, according to the Cauchy--Schwarz inequality, the effective adaptive formulas for the \texttt{DL} parameter given in \cite{HagerZhangCG,DaiKou}, it still has the capability, through an appropriate restriction scheme, to guarantee the sufficient descent condition as well \cite{SBK4OR2}. While the computational performance of $t_k^*$ is quite satisfactory, as demonstrated in \cite{BabaieGhanbariEJOR}---particularly with respect to running time, as expected from its simple analytical form---it is not numerically superior to the parameter $t_k^+$ defined by
\begin{equation*}
    t_k^+ := \dfrac{\|\by_k\|}{\|\bs_k\|} + \dfrac{\bs_k^\top \by_k}{\|\bs_k\|^2},
\end{equation*}
which is obtained by minimizing a sharper upper bound of $\kappa(\bP_{k+1})$ \cite{BabaieGhanbariEJOR}. It can also be regarded as an augmented (larger) version of $t_k^*$.

Another approach is based on directly reducing the distance between $\sigma_k^-$ and $\sigma_k^+$ according to the singular value diversity of $\bP_{k+1}$ in the sense of (\ref{sigma1<1}), it can be shown that
\begin{equation*}
    t_k^\star := \dfrac{\bs_k^\top \by_k}{\|\bs_k\|^2} = \arg\min_{t} \left(\sigma_k^+ - \sigma_k^-\right).
\end{equation*}
It is worth mentioning that, as shown in \cite{BabaieGhanbariOptimization}, $t_k^\star$ represents another optimal choice for the \texttt{DL} parameter, characterized as the minimizer of the distance between the \texttt{DL} search direction and the effective three-term \texttt{CG} direction \texttt{ZZL} \cite{ZhangZhouLiOMS} within a least-squares framework. Here, we address the aforementioned linear reduction of singular value distances, but in terms of the squared versions of the means. Following this idea, it can be observed that
\begin{eqnarray*}
   t_k^\star &=& \arg\min_{t} \left( \textrm{AM}(\mathcal{S}_k)^2 - \textrm{GM}(\mathcal{S}_k)^2 \right) \\
   &=& \arg\min_{t} \left( \textrm{QM}(\mathcal{S}_k)^2 - \textrm{AM}(\mathcal{S}_k)^2 \right) \\
   &=& \arg\min_{t} \left( \textrm{QM}(\mathcal{S}_k)^2 - \textrm{GM}(\mathcal{S}_k)^2 \right).
\end{eqnarray*}

Interestingly, $t_k^+$ can be expressed as the sum of $t_k^\ast$ and $t_k^\star$, while it exhibits higher numerical effectiveness than both. Moreover, computational evidence in the literature confirms the superiority of the formulas derived from the proposed fractional models ($t_k^*$ and $t_k^+$) over $t_k^\star$ obtained from the linear (squared) models \cite{BabaieGhanbariOptimization,BabaieGhanbariEJOR}.

\subsection{Parametric Analysis for the Memoryless Quasi--Newton Algorithms}\label{DL}

As is well known, \texttt{QN} algorithms were originally devised to provide effective schemes for approximating the second-order information of the model (\ref{UP}). In contrast to the Newton method, they have offered significant algorithmic improvements, both theoretically and practically \cite{SunYuan}. In particular, limited-memory \texttt{QN} algorithms have proven remarkably efficient for large-scale instances of (\ref{UP}), successfully addressing the increasing demand for diversity and flexibility in modern optimization tools.

Generally, in the \texttt{QN} algorithms, $\bH_{k+1}$ in (\ref{dk11}) is an approximation of the inverse Hessian, recursively updated by a specific matrix formula \cite{SunYuan}. Hence, well-conditioning of such matrix approximations has always been a significant concern for researchers, which has led to the development of self-scaling \texttt{QN} updating formulas. Among them, the scaled memoryless BFGS (Broyden--Fletcher--Goldfarb--Shanno) and DFP (Davidon--Fletcher--Powell) formulas are the most well-known \cite{OrenSpedicato}, defined by
\begin{eqnarray*}
  \bar{\textbf{Q}}_{k+1} &:=& \theta_k\textbf{I}-\theta_k\dfrac{\bs_k\by_k^\top+\by_k\bs_k^\top}{\bs_k^\top\by_k}+\left(1+\theta_k\dfrac{\by_k^\top\by_k}{\bs_k^\top\by_k}\right)\dfrac{\bs_k\bs_k^\top}{\bs_k^\top\by_k},\\
   \hat{\textbf{Q}}_{k+1}&:=&  \theta_k\textbf{I} - \theta_k\dfrac{\by_k\by_k^\top}{\by_k^\top\by_k}+\dfrac{\bs_k\bs_k^\top}{\bs_k^\top\by_k},
\end{eqnarray*}
starting with $\bar{\textbf{Q}}_{0}=\hat{\textbf{Q}}_{0}=\textbf{I}$. Here, $\theta_k>0$, namely the scaling parameter, is preferably determined by
\begin{equation*}
  \theta_k^* := \dfrac{\bs_k^\top\by_k}{\|\by_k\|^2},
\end{equation*}
for $\bar{\textbf{Q}}_{k+1}$, and by
\begin{equation*}
  \theta_k^\star := \dfrac{\|\bs_k\|^2}{\bs_k^\top\by_k},
\end{equation*}
for $\hat{\textbf{Q}}_{k+1}$ \cite{SBKJOTA3}.

A careful review of the literature reveals the significant potential of scaling parameters to improve the efficiency of memoryless \texttt{QN} algorithms. At the same time, since they serve as scalar estimations of the (inverse) Hessian \cite{BabaieKafaki2015}, their optimal determination remains an open issue requiring further investigation \cite{AndreiOpenProb}. In what follows, we aim to address this problem through the framework of the HM--GM--AM--QM inequality. Accordingly, the following preliminaries are in order.

As discussed in \cite{SBKJOTA3}, the matrix $\bar{\textbf{Q}}_{k+1}$, serving as the SDM of the memoryless BFGS algorithm, possesses $n-2$ eigenvalues equal to $\theta_k$, together with two additional ones, denoted by $\bar{\lambda}_k^-$ and $\bar{\lambda}_k^+$, which are given by
\begin{eqnarray*}
 \bar{\lambda}_k^\pm &:=& \dfrac{1}{2} \left(1+\theta_k\dfrac{\|\by_k\|^2}{\bs_k^\top\by_k}\right)\dfrac{\|\bs_k\|^2}{\bs_k^\top\by_k} \\
  &&\qquad\pm\;\dfrac{1}{2}\sqrt{\left(1+\theta_k\dfrac{\|\by_k\|^2}{\bs_k^\top\by_k}\right)^2\dfrac{\|\bs_k\|^4}{(\bs_k^\top\by_k)^2}-4\theta_k\dfrac{\|\bs_k\|^2}{\bs_k^\top\by_k}},
\end{eqnarray*}
for which, since $\bs_k^\top\by_k>0$, we have
\begin{equation*}\label{lambda<1}
 0<\bar{\lambda}_k^-\leq\theta_k\leq\bar{\lambda}_k^+.
\end{equation*}
Therefore, it is reasonable to focus on reducing the gap between these two limiting eigenvalues to improve the condition number of $\bar{\bQ}_{k+1}$. In this context, by defining $\bar{\mathcal{L}}_k := \{\bar{\lambda}_k^-, \bar{\lambda}_k^+\}$, we can write
\begin{eqnarray*}
   \mathrm{HM}(\bar{\mathcal{L}_k})&=& \dfrac{2\theta_k\left(\bs_k^\top\by_k\right)}{\bs_k^\top\by_k+\theta_k\|\by_k\|^2}, \\
   \mathrm{GM}(\bar{\mathcal{L}_k})&=& \sqrt{\theta_k\dfrac{\|\bs_k\|^2}{\bs_k^\top\by_k}}, \\
   \mathrm{AM}(\bar{\mathcal{L}}_k)&=& \dfrac{1}{2} \left(1+\theta_k\dfrac{\|\by_k\|^2}{\bs_k^\top\by_k}\right)\dfrac{\|\bs_k\|^2}{\bs_k^\top\by_k},\\
   \mathrm{QM}(\bar{\mathcal{L}_k})&=& \sqrt{\dfrac{\|\bs_k\|^2}{\bs_k^\top\by_k}} \times \sqrt{\left(1+\theta_k\dfrac{\|\by_k\|^2}{\bs_k^\top\by_k}\right)^2\dfrac{\|\bs_k\|^2}{\bs_k^\top\by_k}-2\theta_k}.
\end{eqnarray*}
Note that, due to the complexity of the formulaic structure of $\bar{\bQ}_{k+1}$, here we have, to some extent, more complicated formulas for the given means compared with those of $\bP_{k+1}$ as the SDM of the \texttt{DL} algorithm.

Now, in a fractional programming approach, from (\ref{HMGMAMQM2}) we have
\begin{equation*}
    \theta_k^*=\arg\min_{\theta_k} \dfrac{\mathrm{AM}(\bar{\mathcal{L}}_k)}{\mathrm{GM}(\bar{\mathcal{L}}_k)}
    =\arg\min_{\theta_k} \dfrac{\mathrm{GM}(\bar{\mathcal{L}}_k)}{\mathrm{HM}(\bar{\mathcal{L}}_k)}
    =\arg\min_{\theta_k} \dfrac{\mathrm{AM}(\bar{\mathcal{L}}_k)}{\mathrm{HM}(\bar{\mathcal{L}}_k)}.
\end{equation*}
In addition, through a sequence of algebraic manipulations, we can observe that
\begin{equation*}
    \theta_k^*=\arg\min_{\theta_k} \dfrac{\mathrm{QM}(\bar{\mathcal{L}}_k)}{\mathrm{AM}(\bar{\mathcal{L}}_k)}
    =\arg\min_{\theta_k} \dfrac{\mathrm{QM}(\bar{\mathcal{L}}_k)}{\mathrm{GM}(\bar{\mathcal{L}}_k)}
    =\arg\min_{\theta_k} \dfrac{\mathrm{QM}(\bar{\mathcal{L}}_k)}{\mathrm{HM}(\bar{\mathcal{L}}_k)}.
\end{equation*}
It is worth noting that, by linearly reducing the distance between $\bar{\lambda}_k^-$ and $\bar{\lambda}_k^+$ according to the eigenvalue diversity of $\bar{\bQ}_{k+1}$ in the sense of (\ref{lambda<1}), one obtains a multiplicative form of $\theta_k^*$, as shown in \cite{BabaieKafaki2015}, which, however, requires a restriction to ensure its positivity.

Here, we conduct a similar analysis on $\hat{\textbf{Q}}_{k+1}$, the SDM of the memoryless DFP algorithm, which possesses $n-2$ eigenvalues equal to $\theta_k$, along with two additional ones denoted by $\hat{\lambda}_k^-$ and $\hat{\lambda}_k^+$, satisfying
\begin{equation*}\label{lambda<1}
 0 < \hat{\lambda}_k^- \leq \theta_k \leq \hat{\lambda}_k^+.
\end{equation*}
As is well known, there exists a so-called dual relationship between the BFGS and DFP updating formulas \cite{SunYuan}. In this context, by defining $\hat{\mathcal{L}}_k := \{\hat{\lambda}_k^-, \hat{\lambda}_k^+\}$, we can write
\begin{eqnarray*}
   \mathrm{HM}(\hat{\mathcal{L}}_k) &=& \dfrac{2\theta_k\left(\bs_k^\top\by_k\right)^2}{\left(\theta_k\left(\bs_k^\top\by_k\right)+\|\bs_k\|^2\right)\|\by_k\|^2}, \\
   \mathrm{GM}(\hat{\mathcal{L}}_k) &=& \sqrt{\theta_k\dfrac{\bs_k^\top\by_k}{\|\by_k\|^2}}, \\
   \mathrm{AM}(\hat{\mathcal{L}}_k) &=& \dfrac{1}{2}\left(\theta_k+\dfrac{\|\bs_k\|^2}{\bs_k^\top\by_k}\right),\\
   \mathrm{QM}(\hat{\mathcal{L}}_k) &=& \sqrt{  \dfrac{1}{2}\theta_k^2+ \theta_k\dfrac{\|\bs_k\|^2\|\by_k\|^2-\left(\bs_k^\top\by_k\right)^2}{\left(\bs_k^\top\by_k\right)\|\by_k\|^2} + \dfrac{\|\bs_k\|^4}{2\left(\bs_k^\top\by_k\right)^2} }.
\end{eqnarray*}
Hence, in a fractional programming framework, from (\ref{HMGMAMQM2}) it follows that
\begin{equation*}
    \theta_k^\star = \arg\min_{\theta_k}
    \dfrac{\mathrm{AM}(\hat{\mathcal{L}}_k)}{\mathrm{GM}(\hat{\mathcal{L}}_k)}
    = \arg\min_{\theta_k}
    \dfrac{\mathrm{GM}(\hat{\mathcal{L}}_k)}{\mathrm{HM}(\hat{\mathcal{L}}_k)}
    = \arg\min_{\theta_k}
    \dfrac{\mathrm{AM}(\hat{\mathcal{L}}_k)}{\mathrm{HM}(\hat{\mathcal{L}}_k)}.
\end{equation*}
In addition, through a sequence of algebraic manipulations, we can observe that
\begin{equation*}
    \theta_k^\star=\arg\min_{\theta_k} \dfrac{\mathrm{QM}(\hat{\mathcal{L}}_k)}{\mathrm{AM}(\hat{\mathcal{L}}_k)}
    =\arg\min_{\theta_k} \dfrac{\mathrm{QM}(\hat{\mathcal{L}}_k)}{\mathrm{GM}(\hat{\mathcal{L}}_k)}
    =\arg\min_{\theta_k} \dfrac{\mathrm{QM}(\hat{\mathcal{L}}_k)}{\mathrm{HM}(\hat{\mathcal{L}}_k)}.
\end{equation*}

Taken together, further analytical justifications are now provided for several previously proposed parametric formulas of the \texttt{DL} and \texttt{QN} algorithms, reinforcing their theoretical soundness and numerical robustness. The approach is sufficiently flexible to be extended to other variants of the \texttt{DL} algorithm \cite{BabaieKafaki2017}, as well as to the problem of determining the scaling parameter in the general \texttt{QN} updating formulas of the Broyden class \cite{SunYuan}. In particular, the fractional programming models derived from the HM--GM--AM--QM inequality may provide new parametric choices for other memoryless optimization algorithms. In such cases, the simplicity of the analytical structure of the parameter setting strategy is of considerable importance, as it serves as an indicator of efficiency in large-scale models.

\section{Matrix Parametric Configurations for the Dai--Liao Algorithm}\label{MainCG}

A review of the \texttt{CG} literature reveals that the main efforts to enhance the competitiveness of the \texttt{DL} algorithm have primarily focused on developing adaptive strategies for optimally selecting the parameter~$t$ from a specific perspective~\cite{SBKRairo}. Notably, this line of research can be initially attributed to Andrei~\cite{AndreiOpenProb}, who formally introduced it as an open problem. However, all adaptive formulas proposed thus far in the literature have been restricted to scalar settings of~$t$, indicating a noticeable absence of matrix-based adjustments for the \texttt{DL} parameter. As is commonly practiced, the need for diversity in mathematical modeling is often addressed by increasing dimensionality—a classical strategy for broadening theoretical structures. Within the \texttt{DL} framework, several studies have sought to compensate for the loss of accuracy of the algorithm by incorporating higher-order model information~\cite{LotfiHosseini,LuYuanZhan}, while maintaining a reasonable balance between accuracy and computational efficiency.

As a dimensionality expansion approach, we here propose adjusting the \texttt{DL} parameter through a memoryless matrix update, which offers several advantages. First, this modification enhances the flexibility of the \texttt{DL} algorithm while preserving its limited-memory structure. Moreover, it enables potential connections to other \texttt{CG} algorithms, such as the efficient three-term \texttt{ZZL} algorithm~\cite{ZhangZhouLiOMS}. From an accuracy standpoint, the proposed matrix formulation for the \texttt{DL} parameter~$t$ in \eqref{betdl} can be seamlessly incorporated into all variants of the method that employ modified secant equations, which utilize both cost function values and gradient information~\cite{SBKRairo}. In terms of consistency, all \texttt{CG} algorithms derived from the \texttt{DL} framework based on classical \texttt{CG} parameters can also benefit from the proposed matrix adjustment~\cite{SBKRairo}.

To proceed with our scheme, we first reformulate the \texttt{DL} parameter in \eqref{betdl} in matrix form by initially setting $t \longleftarrow t_k\mathbf{I}$, for a given adaptive choice $t_k$. To increase the dimensionality of the parametric adjustment—thereby enhancing the flexibility of the algorithm—we then extend this setting by letting $t \longleftarrow \mathcal{T}_k$, where $\mathcal{T}_k \in \mathbb{R}^{n \times n}$ is a symmetric positive definite matrix, i.e., $\mathcal{T}_k \succ 0$. As a result, the \texttt{DL} search direction can be expressed as
\begin{equation}\label{dkDLT}
     \bd_{0} := -\bg_{0}, \qquad 
     \bd_{k+1} := -\bg_{k+1} + \dfrac{\bg_{k+1}^\top \by_k}{\bd_k^\top \by_k} \bd_k 
     - \dfrac{\bg_{k+1}^\top \bs_k}{\bd_k^\top \by_k} \mathcal{T}_k \bd_k, 
\end{equation}
for all $k \in \mathbb{Z}^+$. Hereafter, we refer to this version of the \texttt{DL} algorithm as the \textbf{Matrix Dai--Liao} (\texttt{MatDL}) method, for which the search direction can also be written as $\bd_{k+1} = -\mathcal{P}_{k+1} \bg_{k+1}$, where
\begin{equation}\label{P}
     \mathcal{P}_{k+1} := \mathbf{I} - \dfrac{\bs_k \by_k^\top}{\bs_k^\top \by_k} 
     + \mathcal{T}_k \dfrac{\bs_k \bs_k^\top}{\bs_k^\top \by_k}, 
     \qquad \forall k \in \mathbb{Z}^+,
\end{equation}
which provides a useful foundation for further matrix-based analysis of the \texttt{DL} method in the context of the proposed matrix-parametric setting.

As a matter of routine, achieving a reasonable balance between accuracy and efficiency is essential when a matrix parameter is embedded into the \texttt{DL} method. Here, to position the \texttt{MatDL} algorithm as an efficient tool that inherits the advantages of the classical \texttt{DL} method, it is crucial to avoid both explicit matrix computations and matrix storage. This can be accomplished by defining $\mathcal{T}_k$ as a rank-one update of the matrix $t_k\mathbf{I}$ using the vector $\by_k$, which explicitly contains recent gradient information and implicitly captures iterative progress. More precisely, we set
\begin{equation}\label{T}
    \mathcal{T}_k := t_k\mathbf{I} + \nu_k\ {\by_k \by_k^\top}, \qquad \forall k \in \mathbb{Z}^+,
\end{equation}
where $t_k$ is the classical \texttt{DL} parameter and $\nu_k$ is an additional nonnegative parameter. It is evident that for $t_k = \nu_k = 0$, \texttt{MatDL} reduces to the \texttt{HS} method, while $\nu_k = 0$ recovers the classical \texttt{DL} method. In what follows, we analytically examine the relationships between $t_k$ and $\nu_k$ to support our primary objectives of flexibility, diversity, and generality. Moreover, to guarantee positive definiteness of $\mathcal{T}_k$, it suffices to impose a positivity assumption on $t_k$.

\subsection{Condition-Oriented Parametric Adjustments for {\tt MatDL}} Since the matrix $\mathcal{T}_k$ defined in \eqref{T} requires only a single vector and two scalar values for its construction, and because the \texttt{MatDL} search direction \eqref{dkDLT} can be computed using only a few inner products, the simplicity criterion is effectively satisfied. At the same time, maintaining an acceptable level of accuracy remains critical. To this end, a straightforward and constructive strategy is to ensure the well-conditioning of the matrix $\mathcal{P}_{k+1}$ defined in \eqref{P}.

As is well known in the matrix computations literature, well-conditioning plays a crucial role in controlling errors along the solution path and enhancing numerical stability \cite{Watkins}. Following the singular value analysis presented in \cite{BabaieGhanbariEJOR}, and noting that under the Wolfe \texttt{LiS} strategy the inequality \eqref{sy} holds, as well as considering the structure of $\mathcal{T}_k$ given in \eqref{T}, the singular values of $\mathcal{P}_{k+1}$ can be given by
\begin{equation*}
    1 \mbox{ with multiplicity } n-2,\ \sigma_k^-,\mbox{ and }\sigma_k^+,
\end{equation*}
in which 
\begin{equation*}
    \sigma_k^-\sigma_k^+ =\dfrac{\|\bs_k\|^2}{\bs_k^\top \by_k}\times\dfrac{\bs_k^\top\mathcal{T}_k\by_k}{\bs_k^\top \by_k},
\end{equation*}
and
\begin{equation*}
    \left(\sigma_k^-\right)^2 + \left(\sigma_k^+\right)^2 = \dfrac{\|\mathcal{T}_k\bs_k\|^2\|\bs_k\|^2}{\left(\bs_k^\top \by_k\right)^2} + \dfrac{\|\bs_k\|^2\|\by_k\|^2}{\left(\bs_k^\top \by_k\right)^2}.
\end{equation*}
Thus,
\begin{equation*}
    \sigma_k^-+\sigma_k^+ =\dfrac{\|\bs_k\|}{\bs_k^\top \by_k}\|\mathcal{T}_k\bs_k+\by_k\|,
\end{equation*}
and consequently,
\begin{equation*}
    \sigma_k^\pm :=\dfrac{\|\bs_k\|}{2\left(\bs_k^\top \by_k\right)}\left(\|\mathcal{T}_k\bs_k+\by_k\|\pm\|\mathcal{T}_k\bs_k-\by_k\|\right).
\end{equation*}

To promote well-conditioning of $\mathcal{P}_{k+1}$, it is essential to reduce its condition number as much as possible, which is classically defined as the ratio of the largest singular value to the smallest singular value of $\mathcal{P}_{k+1}$. Although determining the exact ordering of the singular values of $\mathcal{P}_{k+1}$ typically requires analyzing multiple scenarios—a potentially complex task—it is still possible, as an alternative measure, to control their distribution by minimizing the gap between $\sigma_k^-$ and $\sigma_k^+$. That is,
\begin{equation*}
    \left(t_k^\ast,\nu_k^\ast\right):=\arg\min_{t_k,\nu_k\geq0}\|\mathcal{T}_k\bs_k-\by_k\|^2,
\end{equation*}
for which the following linear system should be solved:
\begin{equation*}
   \left\{\begin{array}{lllll}
 \|\bs_k\|^2t_k &+ & \left(\bs_k^\top \by_k\right)^2\nu_k &=  & \bs_k^\top \by_k,\\
 &&&&\\
  \left(\bs_k^\top \by_k\right)t_k &+ & \left(\bs_k^\top \by_k\right)\|\by_k\|^2\nu_k &= &  \|\by_k\|^2,
\end{array}\right.
\end{equation*}
yielding
\begin{equation}\label{tkthetakstar}
    \left(t_k^\ast,\nu_k^\ast\right) :=\left(0,\dfrac{1}{\bs_k^\top \by_k}\right).
\end{equation}

As an interesting observation, substituting \eqref{tkthetakstar} into \eqref{dkDLT} yields
\begin{equation}\label{dkZZL}
    \nonumber \bd_{k+1} = -\bg_{k+1} + \dfrac{\bg_{k+1}^\top \by_k}{\bd_k^\top \by_k} \bd_k - \dfrac{\bg_{k+1}^\top \bd_k}{\bd_k^\top \by_k} \by_k, \qquad \forall k \in \mathbb{Z}^+,
\end{equation}

which coincides exactly with the effective \texttt{ZZL} direction \cite{ZhangZhouLiOMS}. In other words, through the proposed matrix parameterization of the \texttt{DL} algorithm, the well-known \texttt{ZZL} method can be recovered as a special member of the \texttt{DL} family of \texttt{CG} algorithms. It is worth noting that the \texttt{ZZL} algorithm is not only practically recognized as an efficient method, but also theoretically regarded as an advantageous extension of the \texttt{HS} method. In particular, it satisfies SDC in equality form. Recently, motivated by the merits of the \texttt{ZZL} algorithm, Babaie--Kafaki and Ghanbari proposed an adaptive strategy for setting $t$ in \eqref{betdl} by steering the \texttt{DL} direction as close as possible toward the \texttt{ZZL} direction within the framework of a least-squares model  \cite{SBKRairo}.

In another approach, based on the analysis carried out in \cite{BabaieGhanbariEJOR}, a possible upper bound for the condition number of $\mathcal{P}_{k+1}$ can be expressed as
\begin{equation}\label{upperbound}
\dfrac{ \left(\sigma_k^-\right)^2 + \left(\sigma_k^+\right)^2}{\sigma_k^- \sigma_k^+} = \dfrac{\|\mathcal{T}_k \bs_k\|^2 + \|\by_k\|^2}{\bs_k^\top \mathcal{T}_k \by_k}.
\end{equation}
Noting that
\begin{equation*}
\|\mathcal{T}_k \bs_k\|^2 - 2\left(\bs_k^\top \mathcal{T}_k \by_k\right) + \|\by_k\|^2 = \|\mathcal{T}_k \bs_k - \by_k\|^2 \geq 0,
\end{equation*}
the upper bound in \eqref{upperbound} is minimized when $\mathcal{T}_k \bs_k = \by_k$, which again leads to the parameter setting of \eqref{tkthetakstar}. This observation highlights yet another optimality feature of the parametric scheme \eqref{tkthetakstar} within the \texttt{MatDL} framework defined by \eqref{dkDLT}.

\subsection{Descent-Oriented Parametric Adjustments for {\tt MatDL}} As already mentioned, a primary analytical concern in the context of continuous optimization algorithms is the establishment of the SDC, which is fundamentally linked to their convergence analysis \cite{SugikiNarushimaYabe}. Here, building on the eigenvalue analysis developed in \cite{BabaieGhanbariOMS}, we investigate the descent property of the \texttt{MatDL} algorithm by analyzing the eigenvalues of the following symmetrized version of $\mathcal{P}_{k+1}$:

\begin{eqnarray}
   \nonumber\label{Pb}     \bar{\mathcal{P}}_{k+1} &=& \dfrac{1}{2}\left(\mathcal{P}_{k+1}+\mathcal{P}_{k+1}^\top\right) =\mathbf{I} + \dfrac{\left(\mathcal{T}_k \bs_k - \by_k\right)\bs_k^\top + \bs_k\left(\mathcal{T}_k \bs_k - \by_k\right)^\top}{2\left(\bs_k^\top \by_k\right)}.
\end{eqnarray}
Given the structure of $\mathcal{T}_k$ by \eqref{T}, the eigenvalues of $\bar{\mathcal{P}}_{k+1}$ can be characterized by
\begin{equation*}
 1 \mbox{ with multiplicity } n-2,\ \lambda_k^-,\mbox{ and } \lambda_k^+,
\end{equation*}
where
\begin{equation*}
\lambda_k^- + \lambda_k^+ = 1 + \dfrac{\bs_k^\top \mathcal{T}_k \bs_k}{\bs_k^\top \by_k} > 0,
\end{equation*}
and
\begin{equation}\label{lambda_k-+}
\lambda_k^- \lambda_k^+ = \dfrac{\left(\left(\mathcal{T}_k \bs_k + \by_k\right)^\top\bs_k\right)^2 - \|\mathcal{T}_k \bs_k - \by_k\|^2 \|\bs_k\|^2}{4\left(\bs_k^\top \by_k\right)^2}.
\end{equation}
The positivity of this expression is required to ensure that $\bar{\mathcal{P}}_{k+1} \succ 0$, which, in turn, guarantees the descent condition.

To enforce the positiveness of \eqref{lambda_k-+}, it is necessary to consider manageable forms of the following inequality:
\begin{equation}\label{positiveness}
 \left( \mathcal{T}_k \bs_k + \by_k \right)^\top \bs_k > \| \mathcal{T}_k \bs_k - \by_k \| \| \bs_k \|.
\end{equation}
In this direction, and taking into account the structure of $\mathcal{T}_k$, we observe that
\begin{equation*}
 \left( \mathcal{T}_k \bs_k + \by_k \right)^\top \bs_k > \left( \mathcal{T}_k \bs_k - \by_k \right)^\top \bs_k.
\end{equation*}
Therefore, one strategy to ensure \eqref{positiveness} is to choose the parameters $t_k$ and $\nu_k$ such that
\begin{equation*}
 \left( \mathcal{T}_k \bs_k - \by_k \right)^\top \bs_k  \geq \| \mathcal{T}_k \bs_k - \by_k \|  \| \bs_k \|,
\end{equation*}
which, by the Cauchy--Schwarz inequality, is satisfied only when the vectors $\mathcal{T}_k \bs_k - \by_k$ and $\bs_k$ are parallel. That is,
\begin{equation*}
t_k \bs_k + \left( \nu_k \left( \bs_k^\top \by_k \right) - 1 \right) \by_k \ \boldsymbol{\parallel}\ \bs_k,
\end{equation*}
which again yields $\nu_k^\ast$ defined in \eqref{tkthetakstar}, while $t_k$ can be set as any positive value.


The other approach to guarantee SDC in terms of maintaining \eqref{positiveness} can be founded upon the following inequality:
\begin{equation*}
    \left(\left(\mathcal{T}_k \bs_k + \by_k\right)^\top\bs_k\right)^2>\left(\bs_k^\top\mathcal{T}_k \bs_k\right)^2 \geq t_k^2\|\bs_k\|^4.
\end{equation*}
Hence, the inequality \eqref{positiveness} is ensured when
\begin{equation*}
    t_k^2\|\bs_k\|^2\geq \|\mathcal{T}_k \bs_k - \by_k\|^2.
\end{equation*}
Now, after straightforward algebraic manipulations, it can be seen that 
\begin{equation}\label{tkthetakdescent}
    t_k\geq \dfrac{1-\nu_k\left(\bs_k^\top \by_k\right)}{2} \dfrac{\|\by_k\|^2}{\bs_k^\top \by_k},
\end{equation}
where 
\begin{equation*}
    \nu_k\in\left[0,\nu_k^\ast\right],
\end{equation*}
with $\nu_k^\ast$ defined in \eqref{tkthetakstar}. It is worth highlighting that if $\nu_k = 0$, then the inequality \eqref{tkthetakdescent} reduces to the classical inequality \eqref{DLtau} with $\varsigma = \dfrac{1}{2}$, which is known to guarantee the SDC for the \texttt{DL} algorithm. Conversely, setting $\nu_k = \nu_k^\ast$ (with $t_k = 0$) recovers the \texttt{ZZL} algorithm, which also satisfies SDC.

\subsection{Secant-Oriented Parametric Adjustments for \texttt{MatDL}} It is a well-established practice to incorporate (approximate) higher-order information of an optimization model to enhance accuracy and possibly improve the convergence behavior of the underlying algorithm. Such strategies have been meaningfully developed for the classical \texttt{DL} algorithm, either to increase accuracy by utilizing both cost function values and gradient information, or to establish global convergence without convexity assumptions \cite{SBKRairo}. In this context, secant equations have played a central role, being classically recognized as the foundation of \texttt{QN} algorithms.

Alongside efforts to steer the \texttt{DL} direction toward that of \texttt{QN} methods within a least-squares framework \cite{SBKGhanbari4OR}, other studies have focused on directly equating the \texttt{DL} direction with the \texttt{QN} direction, subsequently determining the \texttt{DL} parameter $t$ in \eqref{betdl} via modified secant equations \cite{LotfiHosseini,LuYuanZhan}. A similar strategy can be applied to the \texttt{MatDL} algorithm with the search direction given in \eqref{dkDLT}. Specifically, we first set
\begin{equation*}
     -\mathbf{B}_{k+1}^{-1}\bg_{k+1} = -\bg_{k+1} + \dfrac{\bg_{k+1}^\top \by_k}{\bd_k^\top \by_k} \bd_k - \dfrac{\bg_{k+1}^\top \bs_k}{\bd_k^\top \by_k} \mathcal{T}_k \bd_k,
\end{equation*}
where the left-hand side represents the general \texttt{QN} direction, with $\mathbf{B}_{k+1} \approx \nabla^2 f(\bx_{k+1})$ as the Hessian approximation. Taking the inner product of both sides with $\bs_k^\top\mathbf{B}_{k+1}$ and performing straightforward algebraic manipulations leads to the following relation:
\begin{equation}\label{TDense}
    \bs_k^\top\mathbf{B}_{k+1}\mathcal{T}_k\bs_k=\dfrac{\bs_k^\top\bg_{k+1}-\bs_k^\top \mathbf{B}_{k+1} \bg_{k+1} +\dfrac{\bg_{k+1}^\top\by_k}{\bs_k^\top \by_k}\left(\bs_k^\top\mathbf{B}_{k+1}\bs_k\right) }
    {\dfrac{\bg_{k+1}^\top\bs_k}{\bs_k^\top \by_k}}.
\end{equation}

Although choosing $\mathbf{B}_{k+1}$ as the exact Hessian or using a dense \texttt{QN} updating formula in \eqref{TDense} can significantly improve accuracy, such choices are computationally expensive, particularly in high-dimensional problems. This challenge can be alleviated by adopting a version of the secant equation \cite{Andreibook}, such as the standard secant condition $\mathbf{B}_{k+1}\bs_k = \by_k$, which enhances the efficiency of the \texttt{MatDL} method through simplified computations. Consequently, considering the rank-one structure of the matrix parameter $\mathcal{T}_k$ defined in \eqref{T}, the following system can be solved to derive adaptive formulas for $t_k$ and $\nu_k$:
\begin{equation}\label{system}
    t_k + \nu_k \|\by_k\|^2 = 1, \quad \text{subject to} \quad t_k, \nu_k \ge 0.
\end{equation}
Alternative forms of this system can be obtained by employing modified secant equations from the literature to simplify \eqref{TDense} \cite{Arazm}.

It is worth noting that the system \eqref{system} can be regarded as a key model in several contexts. First, it can be derived by imposing a conjugacy condition on the \texttt{MatDL} direction \eqref{dkDLT}, specifically that of \eqref{conjugacy0}, which is a well-recognized condition in the literature \cite{Andreibook}. Additionally, the system provides further flexibility for the proposed matrix parametric setting. As illustrative examples, we here suggest the following possible choices consistent with \eqref{system}:
\begin{equation}\label{rhoetta}
    \begin{array}{lcl}
       t_k = \dfrac{\bs_k^\top \by_k}{\varrho \|\bs_k\| \|\by_k\|} & \Longrightarrow & 
       \nu_k = \dfrac{1}{\|\by_k\|^2} \left(1 - \dfrac{\bs_k^\top \by_k}{\varrho \|\bs_k\| \|\by_k\|} \right), \\
       & & \\
       t_k = 1 - \dfrac{\bs_k^\top \by_k}{\eta \|\bs_k\| \|\by_k\|} & \Longrightarrow & 
       \nu_k = \dfrac{\bs_k^\top \by_k}{\eta \|\bs_k\| \|\by_k\|^3},
    \end{array}
\end{equation}

where $\varrho $ and $\eta$ are prespecified constants in $[1, +\infty)$. Moreover, a wide range of solutions for the system \eqref{system} can be constructed analytically, for instance by first selecting $t_k \in [0,1]$ and then setting $\nu_k = \dfrac{1 - t_k}{\|\by_k\|^2}$.

As a final note, to briefly address the global convergence of \texttt{MatDL}, it is important to note that the algorithm possesses a meaningful capability to satisfy the SDC, as discussed in this section. Moreover, since \texttt{MatDL} is rooted in the \texttt{DL} approach, which rarely generates nondescent directions, its convergence can also be established through the analysis carried out in \cite{DaiLiaoNCG}. On the other hand, in the context of Lemma \ref{lemconv}, note that from \eqref{T} we have
\begin{equation*}
    \|\mathcal{T}_k\| \leq t_k + \nu_k\|\by_k\|^2, \qquad \forall k \in \mathbb{Z}^+.
\end{equation*}
So, by imposing boundedness on the parameters $t_k$ and $\nu_k$, the norms of the \texttt{MatDL} search directions can be uniformly bounded as well.


\section{Numerical Experiments}\label{Numerical}

In this section, we investigate the effect of the matrix setting for the
\texttt{DL} parameter, referred to as the \texttt{MatDL} method, and compare it with the classical scalar setting of this parameter. We implement both approaches on the 482 \texttt{CUTEst} test problems of Gould et al.~\cite{CUTEst}, with dimensions ranging from 2 to 9000, and subsequently analyze the obtained results in detail.

The numerical assessment is organized into three successive comparisons,
each corresponding to one of the parametric constructions developed in
Section~\ref{MainCG}. First, we examine the descent-oriented family
obtained by setting
$\nu_k=\texttt{fact}*\,\nu_k^\ast$, where $\nu_k^\ast$ is defined in
\eqref{tkthetakstar} and $\texttt{fact}\in(0,1)$ is a constant, by choosing $t_k$ at the lower bound prescribed
by \eqref{tkthetakdescent}. This comparison evaluates the sensitivity of
\texttt{MatDL} to the parameter \texttt{fact} within a setting designed
to preserve the SDC. Second, the selected
\texttt{fact}-based configuration is compared with the
\texttt{eta}- and \texttt{rho}-based variants obtained from the
secant-oriented system \eqref{system}, as given by \eqref{rhoetta}. These choices are motivated by
the standard secant equation and the successive conjugacy condition
\eqref{conjugacy0}, as discussed in \cite{Andreibook,Arazm}. Third, the
most promising \texttt{MatDL} configurations 
are compared
with the established \texttt{ZZL}, \texttt{DK}, and \texttt{HZ}
methods, while \texttt{DK} and
\texttt{HZ} represent the effective scalar \texttt{DL} parameterizations
associated with the Dai--Kou \cite{DaiKou}, and Hager--Zhang choices \cite{HagerZhangCG}. These three experiments therefore
separate the effects of internal parameter tuning, alternative
theoretically motivated matrix configurations, and comparisons with
efficient \texttt{CG} methods. For all compared methods, we employed the Mor{\'e}--Thuente strong Wolfe
 \texttt{LiS} algorithm \texttt{cvsrch} \cite{cvsrch} with its default
parameter values.

As is widely accepted in the literature, a solver is regarded as more
efficient if it solves a larger proportion of the problems using fewer
function and/or gradient evaluations. To assess this behavior, we employ
the performance profiles of Dolan and Mor\'e \cite{DolM}. For a solver
$s\in\mathfrak{S}$ and a problem $p\in\mathfrak{P}$, where
$\mathfrak{S}$ denotes the set of solvers and $\mathfrak{P}$ the set of
test problems, let $c_{s,p}$ denote the corresponding computational cost.
The performance ratio is defined by
\[
r_{s,p}:=
\frac{c_{s,p}}
{\displaystyle\min_{\mathfrak{s}\in\mathfrak{S}}c_{\mathfrak{s},p}},
\]
and the performance profile of solver $s$ is given by
\[
\rho_s(\tau):=
\frac{1}{|\mathfrak{P}|}
\left|
\left\{
p\in\mathfrak{P}: r_{s,p}\leq\tau
\right\}
\right|.
\]
Thus, $\rho_s(\tau)$ represents the fraction of test problems for which
the cost incurred by solver $s$ is at most $\tau$ times the smallest cost
achieved by any solver in $\mathfrak{S}$. In particular, $\rho_s(1)$
measures efficiency, whereas the limiting value of $\rho_s(\tau)$ as
$\tau\to\infty$ reflects robustness. To simplify the notation, we write
$\rho(\tau):=\rho_s(\tau)$ throughout the remainder of this section.

Note that, when reverse-mode automatic differentiation is employed, the
computational cost of evaluating the gradient is typically only a modest
multiple of that required for evaluating the objective function. For the
\texttt{CUTEst} test set, this multiple is approximately two on average;
see Section~3 of the supplementary material associated with
\texttt{LMBOPT}~\cite{LMBOPT}, available at
\url{https://doi.org/10.5281/zenodo.5607521}. Based on this observation, in the present experiments the two cost measures, \texttt{nf}, denoting
the number of function evaluations, and \texttt{ng}, denoting the number
of gradient evaluations, are combined as
\[
\texttt{nf2g}:=\texttt{nf}+2\,\texttt{ng},
\]
to define a unified measure for comparison within the Dolan--Mor\'e framework  \cite{DolM}.

\begin{figure}[!htbp]
\begin{center}
{\includegraphics[width=15cm]{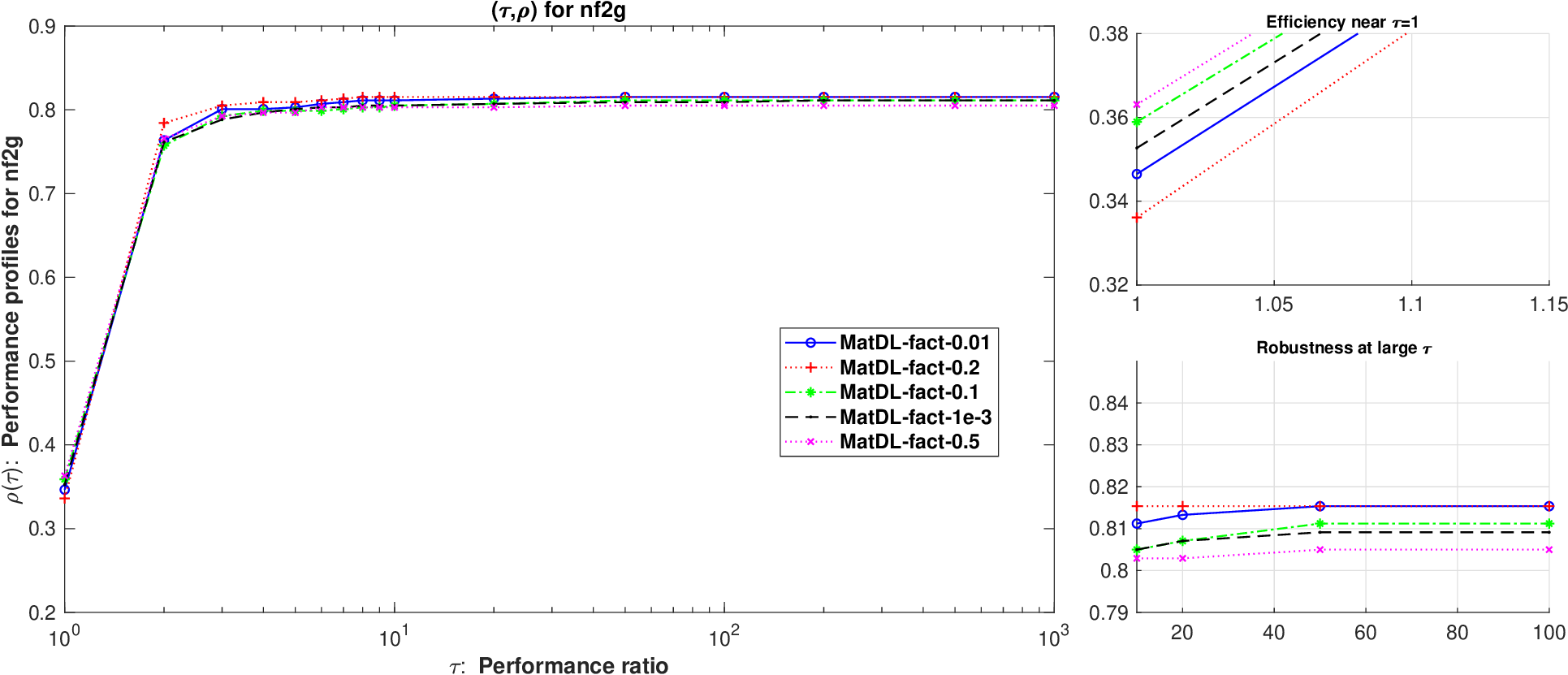}}\\
\end{center}
\caption{Performance profile $\rho(\tau)$ for the \texttt{nf2g} cost measure, independent of an upper bound on the performance ratio $\tau$. The upper-right panel
enlarges the neighborhood of $\tau=1$ to emphasize efficiency, while
the lower-right panel enlarges the large-$\tau$ region to emphasize
robustness. Problems solved by no solver are ignored.}
\label{fig:nf2g-profile1}
\end{figure}

Figure~\ref{fig:nf2g-profile1} shows that no single value of
\texttt{fact} dominates throughout the entire range of performance
ratios. At $\tau=1$, \texttt{MatDL-fact-0.5} attains the largest value
of $\rho(1)$ and is therefore the most efficient variant in terms of
achieving the minimum \texttt{nf2g} cost measure. However, its profile
approaches a lower limiting value as $\tau$ increases, indicating
weaker robustness. In contrast, \texttt{MatDL-fact-0.01} eventually
attains the highest profile and is therefore the most robust variant.
The method \texttt{MatDL-fact-1e-3} also exhibits strong efficiency
near $\tau=1$ while maintaining a high limiting profile. Consequently, both
\texttt{MatDL-fact-0.01} and \texttt{MatDL-fact-1e-3} are retained for
the subsequent comparisons. The former represents the
robustness-oriented choice, whereas the latter mainly provides an
efficiency-oriented alternative while retaining competitive robustness.

\begin{figure}[!htbp]
\begin{center}
{\includegraphics[width=15cm]{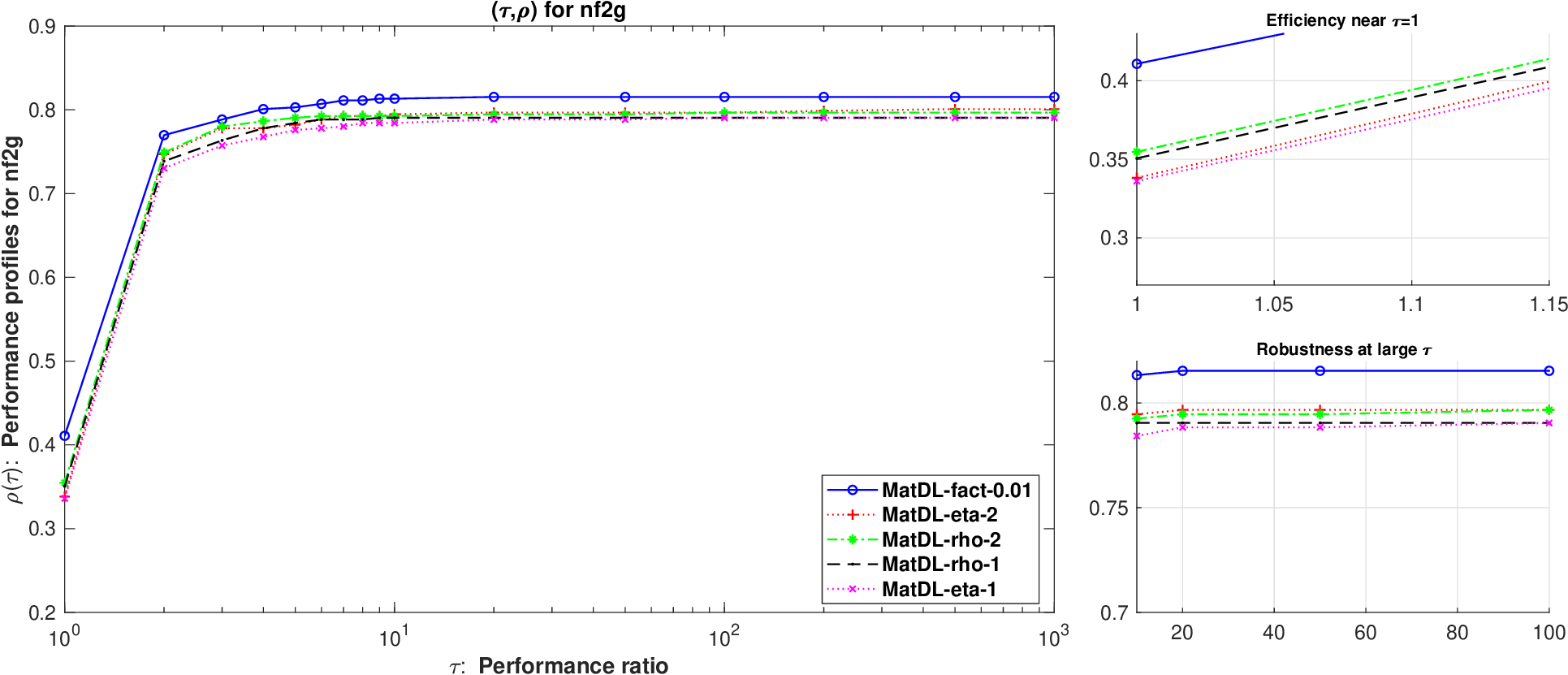}}\\
\end{center}
\caption{Performance profiles $\rho(\tau)$ for the \texttt{nf2g}
cost measure, comparing \texttt{MatDL-fact-0.01} with the
\texttt{eta}- and \texttt{rho}-based variants. The upper-right panel
enlarges the neighborhood of $\tau=1$ to emphasize efficiency, while
the lower-right panel enlarges the large-$\tau$ region to emphasize
robustness. Problems solved by no solver are ignored.}
\label{fig:nf2g-profile3}
\end{figure}

Figure~\ref{fig:nf2g-profile3} compares
\texttt{MatDL-fact-0.01} with the four variants
\texttt{MatDL-eta-1}, \texttt{MatDL-eta-2},
\texttt{MatDL-rho-1}, and \texttt{MatDL-rho-2}.
The enlarged panel near $\tau=1$ shows that
\texttt{MatDL-fact-0.01} attains a substantially larger value of
$\rho(1)$ than all competing variants. It therefore achieves the minimum
\texttt{nf2g}  cost measure on the largest proportion of test problems and is
the most efficient method in this comparison. Among the remaining
variants, \texttt{MatDL-rho-2} exhibits the strongest initial performance,
followed closely by \texttt{MatDL-rho-1}, whereas the two
\texttt{eta}-based variants attain slightly smaller values of
$\rho(1)$. The robustness enlargement for large values of $\tau$ leads
to the same conclusion: \texttt{MatDL-fact-0.01} reaches the highest
limiting profile, while the other four variants level off at noticeably
lower values. In particular, \texttt{MatDL-eta-2} and
\texttt{MatDL-rho-2} are the most robust alternatives, but neither
approaches the limiting profile of \texttt{MatDL-fact-0.01}. Hence,
\texttt{MatDL-fact-0.01} dominates the compared methods in
both efficiency and robustness with respect to the \texttt{nf2g} cost
measure.

\begin{figure}[!htbp]
\centering

\subfloat[
Comparison of \texttt{MatDL-fact-0.01} with
\texttt{ZZL}, \texttt{DK}, and \texttt{HZ}.
\label{fig:nf2g-profile-001}
]{
    \includegraphics[width=0.95\textwidth]{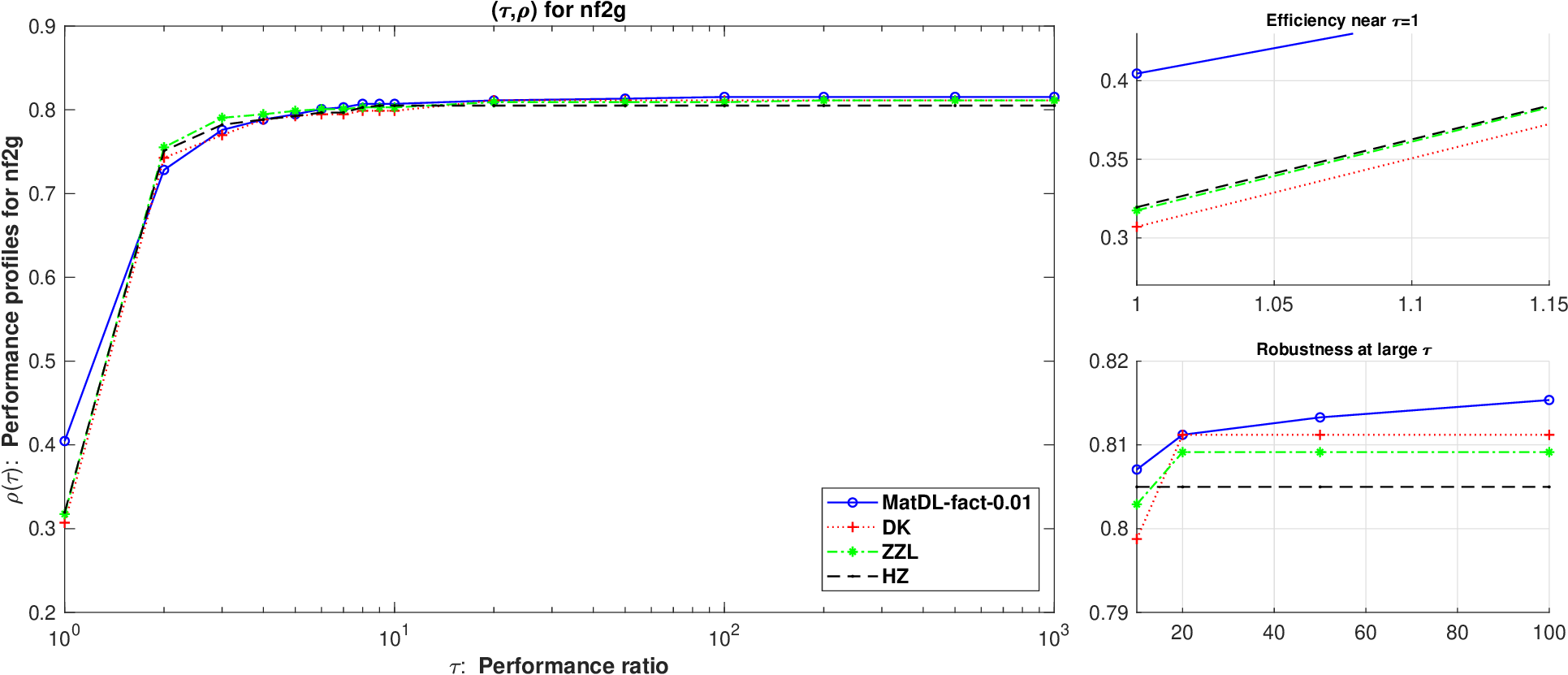}
}

\par\medskip

\subfloat[
Comparison of \texttt{MatDL-fact-1e-3} with
\texttt{ZZL}, \texttt{DK}, and \texttt{HZ}.
\label{fig:nf2g-profile-1e3}
]{
    \includegraphics[width=0.95\textwidth]{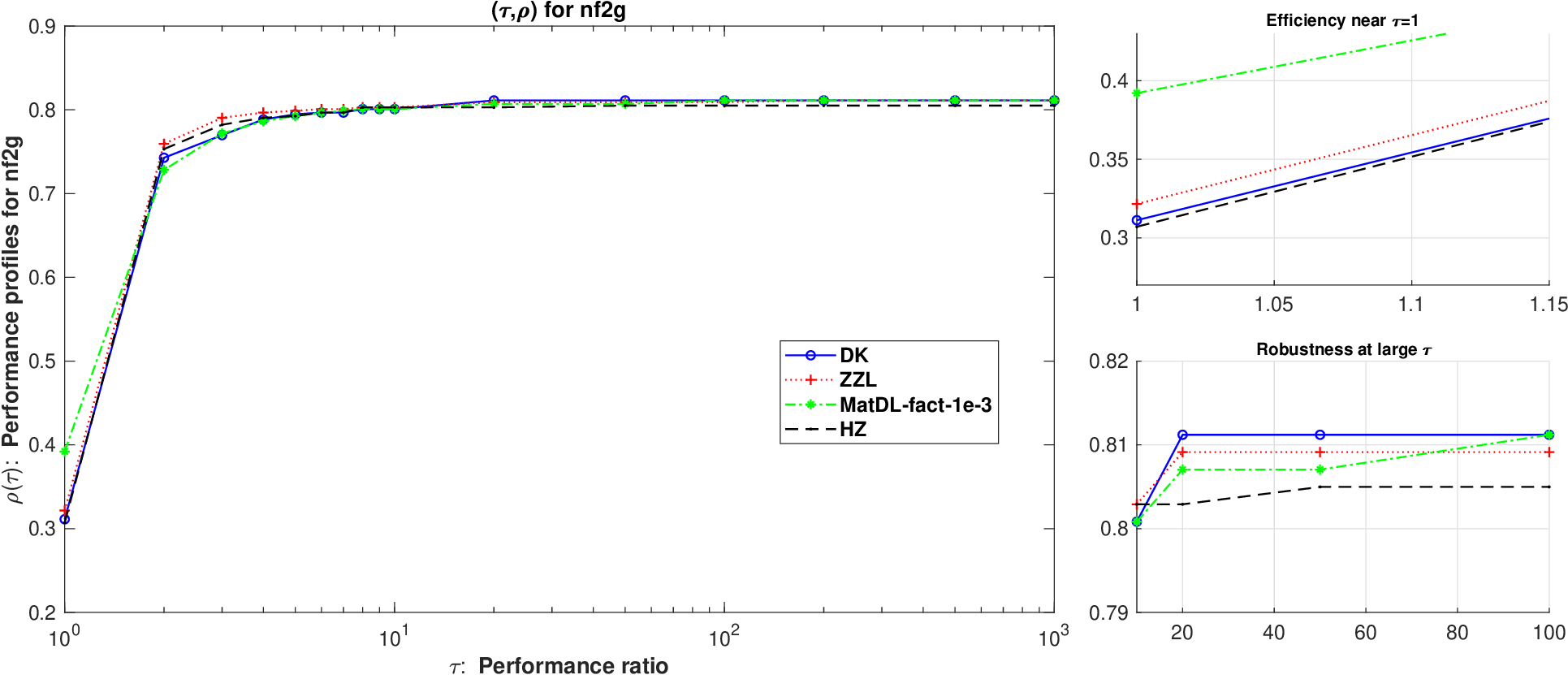}
}

\caption{Performance profiles $\rho(\tau)$ for the \texttt{nf2g}
cost measure. The upper panel presents the robustness-oriented
\texttt{MatDL-fact-0.01} variant, while the lower panel presents the
efficiency-oriented \texttt{MatDL-fact-1e-3} variant. Problems solved
by no solver are ignored.}
\label{fig:nf2g-profile}

\end{figure}

Figure~\ref{fig:nf2g-profile} compares the two selected variants
\texttt{MatDL-fact-0.01} and \texttt{MatDL-fact-1e-3} with
\texttt{ZZL}, \texttt{DK}, and \texttt{HZ}. As shown in
Figure~\ref{fig:nf2g-profile-001}, \texttt{MatDL-fact-0.01} is
particularly competitive for large performance ratios and attains the
highest limiting profile. It therefore represents the
robustness-oriented variant, since it solves the largest proportion of
problems within a finite multiple of the best \texttt{nf2g}  cost measure.
Figure~\ref{fig:nf2g-profile-1e3} shows that
\texttt{MatDL-fact-1e-3} is more competitive near $\tau=1$, indicating
stronger efficiency in terms of achieving, or remaining close to, the
minimum \texttt{nf2g}  cost measure. The profiles of \texttt{ZZL},
\texttt{DK}, and \texttt{HZ} remain competitive over certain ranges of
$\tau$, but none consistently dominates both selected
\texttt{MatDL} variants. Thus, the two parameter choices provide
complementary behavior: \texttt{MatDL-fact-1e-3} is retained as the
efficiency-oriented variant, whereas \texttt{MatDL-fact-0.01} is
retained as the robustness-oriented variant of the numerical comparisons.

\section{Concluding Remarks}\label{Conclusions}

We have investigated parametric settings for the Dai--Liao conjugate
gradient algorithm from both scalar and matrix perspectives. First, the
well-known Harmonic--Geometric--Arithmetic--Quadratic mean inequality
has been employed to provide further optimality characterizations of
previously proposed scalar choices for the Dai--Liao parameter. By
considering the distribution of the singular values of the corresponding
search direction matrix, these choices were interpreted as solutions to
several fractional and difference-based optimization models designed to
promote well-conditioning. A similar analysis has been carried out for the scaling parameters of the scaled memoryless BFGS
and DFP updating formulas, thereby providing a unified justification for the optimality of their classical choices as well.

The main contribution of this study is the extension of the classical
scalar setting of the Dai--Liao parameter to a matrix setting, leading
to the proposed Matrix Dai--Liao (\texttt{MatDL}) method. The rank-one
structure of the matrix parameter preserves the limited-memory
structure of the classical Dai--Liao algorithm and allows the search
direction to be computed using only a small number of vector inner
products, without explicitly forming or storing matrices. Several
parametric configurations have subsequently been developed based on
well-conditioning, the descent property, and the secant equation. The
proposed framework also reveals meaningful connections with existing
efficient conjugate gradient algorithms. In particular, a specific
matrix configuration of the Dai--Liao parameter coincides with an
efficient three-term conjugate gradient direction. The numerical
experiments conducted on 482 unconstrained \texttt{CUTEst} problems,
with dimensions ranging from 2 to 9000, have also provided practical evidence
supporting the proposed matrix parameterization.

The present study also opens several directions for further
investigation. The matrix parameterization may be combined with modified
secant equations incorporating both function and gradient information or
applied to other conjugate gradient methods derived from the
Dai--Liao framework. It would also be worthwhile to investigate richer
low-rank or sparse structures for the matrix parameter. Finally,
extensions to constrained or composite optimization models may further
clarify the practical potential of the proposed matrix-parameterization
approach.

\section*{Declarations}

\noindent
\textbf{Data Availability:} The data that support the findings of this research can be obtained from the corresponding author, subject to a reasonable request.\\

\noindent
\textbf{Conflict of Interest:} No competing interests are declared by the authors.\\

\noindent
\textbf{Ethical Statement:} It is declared that this research did not involve any studies with human participants or animals, and that there are no ethical issues associated with this work.

\textbf{Funding:} Morteza Kimiaei acknowledges financial support from the Austrian Science Foundation under \url{https://doi.org/10.55776/PAT2747625}.




\begin{thebibliography}{10}

\bibitem{AndreiComparison}
Neculai Andrei.
\newblock Numerical comparison of conjugate gradient algorithms for
  unconstrained optimization.
\newblock {\em Studies in Informatics and Control}, 16(4):333{--}352, 2007.

\bibitem{AndreiOpenProb}
Neculai Andrei.
\newblock Open problems in conjugate gradient algorithms for unconstrained
  optimization.
\newblock {\em Bulletin of the Malaysian Mathematical Sciences Society},
  34(2):319{--}330, 2011.

\bibitem{Andreibook}
Neculai Andrei.
\newblock {\em Modern Numerical Nonlinear Optimization}.
\newblock Springer, Cham, 2022.

\bibitem{Arazm}
Mohammad~Reza Arazm, Saman Babaie{--}Kafaki, and Reza Ghanbari.
\newblock An extended {Dai--Liao} conjugate gradient method with global
  convergence for nonconvex functions.
\newblock {\em Glasnik Matematicki}, 52(2):361--375, 2017.

\bibitem{SBK4OR2}
Saman Babaie{--}Kafaki.
\newblock On the sufficient descent condition of the {H}ager{--}{Z}hang
  conjugate gradient methods.
\newblock {\em 4OR}, 12(3):285--292, 2014.

\bibitem{BabaieKafaki2015}
Saman Babaie{--}Kafaki.
\newblock A modified scaling parameter for the memoryless {BFGS} updating
  formula.
\newblock {\em Numerical Algorithms}, 72(2):425--433, 2015.

\bibitem{SBKJOTA3}
Saman Babaie{--}Kafaki.
\newblock On optimality of the parameters of self-scaling memoryless
  quasi--{N}ewton updating formulae.
\newblock {\em Journal of Optimization Theory and Applications},
  167(1):91--101, 2015.

\bibitem{SBKRairo}
Saman Babaie{--}Kafaki.
\newblock A survey on the {Dai--Liao} family of nonlinear conjugate gradient
  methods.
\newblock {\em RAIRO--Operations Research}, 57:43--58, 2023.

\bibitem{BabaieGhanbariOptimization}
Saman Babaie{--}Kafaki and Reza Ghanbari.
\newblock Two optimal {Dai--Liao} conjugate gradient methods.
\newblock {\em Optimization}, 64(11):2277--2287, 2015.

\bibitem{SBKGhanbari4OR}
Saman Babaie{--}Kafaki and Reza Ghanbari.
\newblock A class of adaptive {Dai--Liao} conjugate gradient methods based on
  the scaled memoryless {BFGS} update.
\newblock {\em 4OR}, 15(1):85--92, 2017.

\bibitem{BabaieKafaki2017}
Saman Babaie{--}Kafaki and Reza Ghanbari.
\newblock An optimal extension of the {Polak--Ribi\`{e}re--Polyak} conjugate
  gradient method.
\newblock {\em Numerical Functional Analysis and Optimization},
  38(9):1115--1124, 2017.

\bibitem{BabaieGhanbariOMS}
Saman Babaie–Kafaki and Reza Ghanbari.
\newblock A descent family of {Dai–Liao} conjugate gradient methods.
\newblock {\em Optimization Methods and Software}, 29(3):583–591, 2013.

\bibitem{BabaieGhanbariEJOR}
Saman Babaie–Kafaki and Reza Ghanbari.
\newblock The {Dai–Liao} nonlinear conjugate gradient method with optimal
  parameter choices.
\newblock {\em European Journal of Operational Research}, 234(3):625–630,
  2014.

\bibitem{DaiKou}
Yu{--}Hong Dai and Chang{--}Xing Kou.
\newblock A nonlinear conjugate gradient algorithm with an optimal property and
  an improved wolfe line search.
\newblock {\em SIAM Journal on Optimization}, 23(1):296--320, 2013.

\bibitem{CGConvDai}
Yuhong Dai, Jiye Han, Guanghui Liu, Defeng Sun, Hongxia Yin, and Ya–Xiang
  Yuan.
\newblock Convergence properties of nonlinear conjugate gradient methods.
\newblock {\em SIAM Journal on Optimization}, 10(2):345–358, 2000.

\bibitem{DaiLiaoNCG}
Yu–Hong Dai and Li–Zhi Liao.
\newblock New conjugacy conditions and related nonlinear conjugate gradient
  methods.
\newblock {\em Applied Mathematics and Optimization}, 43(1):87{--}101, 2001.

\bibitem{DargahiMMA}
Fatemeh Dargahi, Saman Babaie{--}Kafaki, Zohre Aminifard, and Mohammad
  Hajian{--}Berenjestanaki.
\newblock Solving an augmented nonnegative matrix factorization model by
  modified scaled nonmonotone memoryless {BFGS} methods devised based on the
  ellipsoid vector norm.
\newblock {\em Mathematical Methods in the Applied Sciences},
  48(8):9088{--}9097, 2025.

\bibitem{DennisWolkowicz}
Jr. Dennis, John~E. and Henry Wolkowicz.
\newblock Sizing and least-change secant methods.
\newblock {\em SIAM Journal on Numerical Analysis}, 30(5):1291--1314, 1993.

\bibitem{DolM}
Elizabeth~D. Dolan and Jorge~J. Moré.
\newblock Benchmarking optimization software with performance profiles.
\newblock {\em Mathematical Programming}, 91(2):201--213, jan 2002.

\bibitem{Elden}
Lars Eld{\'e}n.
\newblock {\em Matrix Methods in Data Mining and Pattern Recognition}.
\newblock SIAM, Philadelphia, PA, 2007.

\bibitem{GN}
Jean~Charles Gilbert and Jorge Nocedal.
\newblock Global convergence properties of conjugate gradient methods for
  optimization.
\newblock {\em SIAM Journal on Optimization}, 2(1):21–42, 1992.

\bibitem{CUTEst}
Nicholas I.~M. Gould, Dominique Orban, and Philippe~L. Toint.
\newblock {CUTEst}: a constrained and unconstrained testing environment with
  safe threads for mathematical optimization.
\newblock {\em Computational Optimization and Applications}, 60:545--557, 2015.

\bibitem{HagerZhangCG}
William~W. Hager and Hongchao Zhang.
\newblock A new conjugate gradient method with guaranteed descent and an
  efficient line search.
\newblock {\em SIAM Journal on Optimization}, 16(1):170–192, 2005.

\bibitem{HagerZhang}
William~W. Hager and Hongchao Zhang.
\newblock A survey of nonlinear conjugate gradient methods.
\newblock {\em Pacific Journal of Optimization}, 2(1):35--58, 2006.

\bibitem{HS}
Magnus~R. Hestenes and Eduard Stiefel.
\newblock Methods of conjugate gradients for solving linear systems.
\newblock {\em Journal of Research of the National Bureau of Standards},
  49(6):409--436, 1952.

\bibitem{Kimiaei2026}
Morteza Kimiaei, Saman Babaie–Kafaki, Vyacheslav Kungurtsev, and Mahsa
  Yousefi.
\newblock Restarting two-term nonlinear conjugate gradient methods based on a
  finite precision arithmetic analysis.
\newblock {\em Journal of Computational and Applied Mathematics}, 488:117856,
  2026.

\bibitem{LMBOPT}
Morteza Kimiaei, Arnold Neumaier, and Behzad Azmi.
\newblock {LMBOPT}: a limited memory method for bound-constrained optimization.
\newblock {\em Mathematical Programming Computation}, 14(2):271--318, 2022.

\bibitem{LotfiHosseini}
Mina Lotfi and S.~Mohammad Hosseini.
\newblock An efficient {Dai--Liao} type conjugate gradient method by
  reformulating the {CG} parameter in the search direction equation.
\newblock {\em Journal of Computational and Applied Mathematics}, 371:112708,
  2020.

\bibitem{LuYuanZhan}
Junyu Lu, Gonglin Yuan, and Zhan Wang.
\newblock A modified {Dai--Liao} conjugate gradient method for solving
  unconstrained optimization and image restoration problems.
\newblock {\em Journal of Applied Mathematics and Computing}, 68(2):681–703,
  2021.

\bibitem{cvsrch}
Jorge~J. Mor{\'e} and David~J. Thuente.
\newblock Line search algorithms with guaranteed sufficient decrease.
\newblock {\em ACM Transactions on Mathematical Software (TOMS)},
  20(3):286--307, 1994.

\bibitem{OrenSpedicato}
Saul~S. Oren and Emilio Spedicato.
\newblock Optimal conditioning of self-scaling variable metric algorithms.
\newblock {\em Mathematical Programming}, 10(1):70--90, 1976.

\bibitem{Peressini}
Anthony~L. Peressini, Francis~E. Sullivan, and Jr. Uhl, John~J.
\newblock {\em The Mathematics of Nonlinear Programming}.
\newblock Springer, New York, 1988.

\bibitem{Arashi}
A.~K. Md.~Ehsanes Saleh, Mohammad Arashi, and B.~M.~Golam Kibria.
\newblock {\em Theory of Ridge Regression Estimation with Applications}.
\newblock John Wiley \& Sons, Hoboken, NJ, 2019.

\bibitem{Inequalities}
Hayk Sedrakyan and Nairi Sedrakyan.
\newblock {\em Algebraic Inequalities}.
\newblock Springer, Cham, 2018.

\bibitem{SugikiNarushimaYabe}
Kaori Sugiki, Yasushi Narushima, and Hiroshi Yabe.
\newblock Globally convergent three-term conjugate gradient methods that use
  secant conditions and generate descent search directions for unconstrained
  optimization.
\newblock {\em Journal of Optimization Theory and Applications},
  153(3):733–757, 2011.

\bibitem{SunYuan}
Wenyu Sun and Ya{--}Xiang Yuan.
\newblock {\em Optimization Theory and Methods: Nonlinear Programming}.
\newblock Springer, New York, 2006.

\bibitem{Watkins}
David~S. Watkins.
\newblock {\em Fundamentals of Matrix Computations}.
\newblock John Wiley \& Sons, New York, 2002.

\bibitem{ZhangZhouLiOMS}
Li~Zhang, Weijun Zhou, and Donghui Li.
\newblock Some descent three-term conjugate gradient methods and their global
  convergence.
\newblock {\em Optimization Methods and Software}, 22(4):697–711, 2007.

\end{thebibliography}
\end{document}